%% file: GP-20211210a.tex
\documentclass[leqno,12pt]{article}
\usepackage{amsmath,amssymb,rotating,a4wide,color}
\usepackage{wasysym}
\usepackage{hyperref}

\usepackage{tikz}
\usetikzlibrary{circuits.logic.US,circuits.logic.IEC,fit}
\usetikzlibrary{backgrounds}
\usetikzlibrary{cd}
\usepackage{subcaption}
\usepackage{lscape}
\usepackage{pgfplots}

\def\FigZ{\Placeholder}
\def\FigA{\Placeholder}
\def\FigB{\Placeholder}
\def\FigC{\Placeholder}
\def\FigD{\Placeholder}
\input{GP-Figures.tex}

\usepackage{array}
\usepackage[all,cmtip]{xy}

\usepackage{verbatim}

\newdir^{ (}{{}*!/-3pt/\dir^{(}}
\newdir_{ (}{{}*!/-3pt/\dir_{(}}

\entrymodifiers={+!!<0pt,\fontdimen22\textfont2>}

\def\bigtimes{\mathop{\raise-2pt\hbox{\huge$\times$}}}

\newbox\circbulletbox
\setbox\circbulletbox=\hbox{\,{\large$\circledcirc$}\kern-8.7pt\raise .6pt\hbox{$\bullet$}\,}

\let\le\leqslant
\let\ge\geqslant

\def\circVbig{\hbox{\text{\it\r{V}}}}
\def\circVscript{\hbox{\scriptsize\text{\it\r{V}}}}
\def\circVscriptscript{\mbox{\tiny\text{\it\r{V}}}}

\def\circVlimits_#1^#2{{\mathchoice%
   {\circVbig{}^{\kern2pt #2}_{\kern-2pt #1}}%
   {\circVbig{}^{\kern2pt #2}_{\kern-2pt #1}}%
   {\scriptstyle\circVscript{}^{\kern1.7pt #2}_{\kern-1pt #1}}%
   {\scriptscriptstyle\circVscriptscript{}^{\kern1.5pt #2}_{\kern-1pt #1}}%
   }}
\def\circVr_#1{\circVlimits_#1^r}
\def\circVs_#1{\circVlimits_#1^s}

\def\circWbig{\hbox{\text{\it\r{W}}}}
\def\circWscript{\hbox{\scriptsize\text{\it\r{W}}}}
\def\circWscriptscript{\mbox{\tiny\text{\it\r{W}}}}

\def\circWlimits_#1^#2{{\mathchoice%
   {\circWbig{}^{\kern2pt #2}_{\kern-2pt #1}}%
   {\circWbig{}^{\kern2pt #2}_{\kern-2pt #1}}%
   {\scriptstyle\circWscript{}^{\kern1.7pt #2}_{\kern-1pt #1}}%
   {\scriptscriptstyle\circWscriptscript{}^{\kern1.5pt #2}_{\kern-1pt #1}}%
   }}

\def\OM{\mathchoice
  {\rlap{\kern3.2pt$\overline{\phantom{L}}$}M}
  {\rlap{\kern3.2pt$\overline{\phantom{L}}$}M}
  {\rlap{\kern2.4pt$\scriptstyle\overline{\phantom{L}}$}M}
  {\rlap{\kern1.8pt$\scriptscriptstyle\overline{\phantom{L}}$}M}}

\def\mycirc{{\kern1pt\circ\kern2pt}}
\def\charact{\mathop{\rm char}\nolimits}

\def\sm{{\rm sm}}

\def\Spec{\mathop{\rm Spec}\nolimits}

\def\red{{\rm red}}

\let\phi\varphi
\let\theta\vartheta
\let\epsilon\varepsilon
\let\setminus\smallsetminus

\let\emptyset\varnothing

\newtheorem{Thm}{Theorem}[section]
\newtheorem{Prop}[Thm]{Proposition}
\newtheorem{Lem}[Thm]{Lemma}

\newtheorem{Def}[Thm]{Definition}

\newtheorem{Rem}[Thm]{Remark}

\numberwithin{Thm}{section}
\def\UseTheoremCounterForNextEquation{\setcounter{equation}{\value{Thm}}\addtocounter{Thm}{1}}

\def\qed{{\hskip0pt\unskip\unskip\nobreak\hfil\penalty50
          \hskip1em\hbox{}\nobreak\hfil
           {$\square$}
          \parfillskip=0pt\finalhyphendemerits=0
          \par}\medskip}

\newenvironment{Proof}
               {\noindent{\bf Proof.}\ }
               {\qed}

               {\noindent{\bf Proof of #1.}\ }
               {\qed}

\newcommand{\BA}{{\mathbb{A}}}

\newcommand{\BP}{{\mathbb{P}}}

\newcommand{\BZ}{{\mathbb{Z}}}

\newcommand{\Fm}{{\mathfrak{m}}}

\newcommand{\CC}{{\cal C}}
\newcommand{\CD}{{\cal D}}

\newcommand{\CL}{{\cal L}}

\newcommand{\CO}{{\cal O}}
\newcommand{\CP}{{\cal P}}
\newcommand{\CQ}{{\cal Q}}
\newcommand{\CR}{{\cal R}}

\newcommand{\CU}{{\cal U}}
\newcommand{\CV}{{\cal V}}

\newcommand{\CX}{{\cal X}}

\newbox\mybox
\def\arrover#1{\mathrel{
       \setbox\mybox=\hbox spread 1.4em
              {\hfil$\scriptstyle#1$\hfil}
       \vbox{\offinterlineskip\copy\mybox
             \hbox to\wd\mybox{\rightarrowfill}}}}

\def\larrover#1{\mathrel{
       \setbox\mybox=\hbox spread 1.4em
              {\hfil$\scriptstyle#1\vphantom{g}$\hfil}
       \vbox{\offinterlineskip\copy\mybox
             \hbox to\wd\mybox{\leftarrowfill}}}}

\def\ontoover#1{\mathrel{
       \setbox\mybox=\hbox spread 1.4em
              {\hfil$\scriptstyle#1\vphantom{g}$\hfil}
       \vbox{\offinterlineskip\copy\mybox
             \hbox to\wd\mybox{\rightarrowfill\hskip-2.8mm
                               $\rightarrow$}}}}
\def\leftontoover#1{\mathrel{
       \setbox\mybox=\hbox spread 1.4em
              {\hfil$\scriptstyle#1\vphantom{g}$\hfil}
       \vbox{\offinterlineskip\copy\mybox
             \hbox to\wd\mybox{$\leftarrow$\hskip-2.8mm
                               \leftarrowfill}}}}
\let\longto\longrightarrow
\let\into\hookrightarrow
\let\onto\twoheadrightarrow
\def\longonto{\ontoover{\ }}
\def\longinto{\lhook\joinrel\longrightarrow}

\def\isoto{\mathrel{
       \setbox\mybox=\hbox spread 0.9em
              {\hfil$\scriptstyle\sim$\hfil}
       \vbox{\offinterlineskip\copy\mybox
             \hbox to\wd\mybox{\rightarrowfill}}}}

\begin{document}

\title{\strut
\vskip-80pt
Reduction of Hyperelliptic Curves in Characteristic $\not=2$}

\author{
\begin{minipage}{.3\hsize}
Tim Gehrunger\\[12pt]
\small Department of Mathematics \\
ETH Z\"urich\\
8092 Z\"urich\\
Switzerland \\
tim.gehrunger@math.ethz.ch\\[9pt]
\end{minipage}
\qquad
\begin{minipage}{.3\hsize}
Richard Pink\\[12pt]
\small Department of Mathematics \\
ETH Z\"urich\\
8092 Z\"urich\\
Switzerland \\
pink@math.ethz.ch\\[9pt]
\end{minipage}
}

\date{\today}

\maketitle

\begin{abstract}
Let $K$ be the quotient field of a discrete valuation ring $R$ with residue characteristic $\not=2$, and let $C$ be a hyperelliptic curve over~$K$. We assume that all geometric branch points of the double covering $C\onto\BP^1_K$ are rational and mark both $C$ and $\BP^1_K$ with these branch points. After possibly replacing $R$ by a ramified extension of degree~$2$, we give a direct construction for the stable model of $C$ as a marked curve over~$R$. We deduce that the closed fiber of this stable model is determined completely by the closed fiber of the stable model of the marked~$\BP^1_K$. In particular, the dual graph and other information for the former can be read off directly from the corresponding information for the latter.
\end{abstract}

{\renewcommand{\thefootnote}{}
}

\tableofcontents
\newpage


\section{Introduction}
\label{Intro}

\hskip\parindent{\bf Setup:} 
Let $R$ be a discrete valuation ring with quotient field $K$ and residue field $k$ of characteristic $\neq 2$. Let $C$ be a hyperelliptic curve over~$K$, that is, a double covering $\pi\colon C \onto \bar C$ of a rational curve~$\bar C$.
Since $K$ has characteristic $\not=2$, this covering is only tamely ramified. By the Hurwitz formula it is therefore ramified at exactly $2g+2$ geometric points, where $g$ is the genus of~$C$. After replacing $K$ by a finite extension, if necessary, we assume that these geometric points are all defined over~$K$.
Let $P_1,\ldots,P_{2g+2} \in C(K)$ be these points, and let $\bar P_1,\ldots,\bar P_{2g+2} \in \bar C(K)$ denote their images under~$\pi$.

By the theory of stable marked curves of genus $0$ there exists a stable marked curve $(\bar\CC,\bar\CP_1,\ldots,\bar\CP_{2g+2})$ over $S:=\Spec R$ with generic fiber $(\bar C,\bar P_1,\ldots,\bar P_{2g+2})$, which is unique up to unique isomorphism and can be computed by a simple explicit algorithm. 
The theory of stable marked curves of arbitrary genus says that, after possibly replacing $R$ by a finite extension, there exists a stable marked curve $(\CC,\CP_1,\ldots,\CP_{2g+2})$ over $S$ with generic fiber $(C,P_1,\ldots,P_{2g+2})$, and that it is unique up to unique isomorphism, once it exists. 

\medskip
{\bf Main Results:} 
In this paper we give a direct construction of $(\CC,\CP_1,\ldots,\CP_{2g+2})$ from $(\bar\CC,\bar\CP_1,\ldots,\bar\CP_{2g+2})$ and the given covering $\pi\colon C\onto\bar C$ over~$K$. In fact, we realize $\CC$ as a double covering $\Pi\colon {\CC\onto\bar\CC}$ which extends~$\pi$, using explicit local computations everywhere on~$\bar\CC$. The construction shows that, once all branch points of $\pi$ are defined over~$K$, an additional ramified extension of $R$ of degree at most $2$ is enough for the existence of $(\CC,\CP_1,\ldots,\CP_{2g+2})$.

The method also shows that the special fiber $(C_0,P_{0,1},\ldots,P_{0,2g+2})$ of $(\CC,\CP_1,\ldots,\CP_{2g+2})$ is determined completely by the special fiber $(\bar C_0,\bar P_{0,1},\ldots,\bar P_{0,2g+2})$ of $(\bar\CC,\bar\CP_1,\ldots,\bar\CP_{2g+2})$. 
In particular, every irreducible component $X$ of $C_0$ is either a double covering of an irreducible component $\bar X$ of~$\bar C_0$ with specified branch points, or maps isomorphically to~$\bar X$. Specifically $X$ is smooth without self-intersection, and of genus $\le1$ or hyperelliptic.
Moreover, one can deduce the dual graph of $C_0$ and the genus of each irreducible component directly from the dual graph of $\bar C_0$ and the location of its marked points.
{}From this information one can also read off the reduction behavior of the Jacobian variety of~$C$.

%
%
%

\medskip
{\bf Explanation:} 
For a different perspective on these results, observe that by the uniqueness of the stable model, the hyperelliptic involution of $C$ over~$\bar C$ extends uniquely to an involution $\sigma$ of~$\bar\CC$. As a byproduct, our construction yields a natural isomorphism $\CC/\langle\sigma\rangle \cong \bar\CC$.
While it is known in much greater generality that $\CC/\langle\sigma\rangle$ is a semistable curve, 
the fact that this quotient with its induced marked points is already 
stable and hence isomorphic to $(\bar\CC,\bar\CP_1,\ldots,\bar\CP_{2g+2})$ is the main reason why the simple construction in this paper is possible. 

This fact depends crucially on the assumption $\charact(k)\not=2$. 
To see how, we determine when $C$ has good reduction.
Assume first that the branch points $\bar P_1,\ldots,\bar P_{2g+2}$ extend to disjoint sections of a smooth model $\bar\CC$ of~$\bar C$. In terms of a suitable coordinate $x$ on $\bar\CC\cong\BP^1_S$, the curve $C$ is given by the equation $y^2=\prod_i(x-\xi_i)$ with disjoint sections $\xi_1,\ldots,\xi_n$. Since $2$ is invertible over~$S$, the same equation then directly defines a smooth hyperelliptic curve $\CC$ over~$S$ with $\CC/\langle\sigma\rangle \cong\bar\CC$, showing that the unmarked curve $C$ has good reduction. 
Conversely assume that the unmarked curve $C$ has good reduction with smooth model~$\CC$. Then $\bar\CC := \CC/\langle\sigma\rangle$ is a smooth model of~$\bar C$ and, since $2$ is invertible over~$S$, the quotient morphism $\CC\onto\bar\CC$ is only tamely ramified in every fiber. By the Hurwitz formula again, the ramification locus therefore consists of exactly $2g+2$ geometric points in every fiber. But by the purity of the branch locus, all the $2g+2$ branch points in the special fiber are reductions of the $2g+2$ branch points in the generic fiber; hence the latter extend to disjoint sections of~$\bar\CC$. 
Together this shows that $C$ has good reduction if and only the stable marked curve $(\bar C,\bar P_1,\ldots,\bar P_{2g+2})$ has good reduction, and then $\CC$ is a double covering of~$\bar\CC$.

\medskip
{\bf Contrast:} 
In the case $\charact(k)=2$ the situation is very different. If $K$ has characteristic~$0$, one can still define the stable marked curves $(\bar\CC,\bar\CP_1,\ldots,\bar\CP_{2g+2})$ and $(\CC,\CP_1,\ldots,\CP_{2g+2})$, and since the quotient of the latter by the hyperelliptic involution is a semistable marked curve, the double covering $\pi$ extends uniquely to a morphism $\CC\onto\bar\CC$. But this morphism cannot be finite if $g\ge1$. The reason for this is that a hyperelliptic curve of genus $g$ in characteristic $2$ is wildly ramified at every branch point; hence the number of branch points is at most $g+1$. Thus even if $C$ has good reduction, the branch points in the generic fiber must come together at least pairwise in the special fiber of a smooth model, and separating them in the stable marked model $(\CC,\CP_1,\ldots,\CP_{2g+2})$ requires blowups. 
The situation can become much more tricky if $C$ has bad reduction.

If $\charact(k)=2$, the relation between $(\bar\CC,\bar\CP_1,\ldots,\bar\CP_{2g+2})$ and $(\CC,\CP_1,\ldots,\CP_{2g+2})$ is therefore quite complicated. The case that the stable marked curve $(\bar C,\bar P_1,\ldots,\bar P_{2g+2})$ has good reduction is described by Lehr and Matignon \cite{LehrMatignon2006}, while the general case is still open. The special case $g=1$ was treated in the master thesis \cite{Gehrunger2020} of the first author, on which the present article is partly based. We hope to address the general case in a future paper. 

\medskip
{\bf Relation with other work:} 
Semistable models of hyperelliptic curves in residue characteristic $\not=2$ were first constructed by Bosch \cite{Bosch1980} using the language of rigid geometry, who also described the structure of the special fiber to some extent. 
For an English summary of this see Fresnel and van der Put \cite[Examples 5.5.4]{FresnelvanderPut2004}. 
%
Next, Mochizuki \cite[\S3.13]{Mochizuki1995}, respectively Bouw and Wewers \cite[\S4]{BouwWewers2017}, repeat the same construction in the language of algebraic geometry and for more general finite coverings of curves.
%
%
All of these end up with the same model $\CC$ as we do, though without discussing the stability of the marked curve $(\CC,\CP_1,\ldots,\CP_{2g+2})$.

Kausz \cite[\S4]{Kausz1999} and Srinivasan \cite[\S2]{Srinivasan2015} again use algebraic geometry, but concentrate on \emph{regular} models of hyperelliptic curves. This requires blowing up $\bar\CC$, and each paper shows that the normalization of such a regular blowup is a regular semistable model of~$C$. They also discuss the structure of the special fiber to some extent, but not questions of stability.
%
%
%

The paper \cite{DDMM} by Dokchitser, Dokchitser, Maistret, and Morgan 
is closest in scope to ours, but came to our attention only after we had finished our work. They construct more or less the same covering $\CC\onto\bar\CC$ as we do and obtain comparable results about its special fiber --- with many more details besides --- but replace the well-established language of stable marked curves by their notion of ``cluster pictures''.

Thus, most of what we do has already been done in one form or other. But we feel that it will be of benefit to publish an independent treatment based on stable marked curves.

In arbitrary residue characteristic, the stable reduction of an unmarked hyperelliptic curve of genus $2$ was described completely in Liu \cite[Thm.\,1]{Liu1993}.
A special case of the reduction of hyperelliptic curves to residue characteristic $2$ is treated in Lehr and Matignon \cite{LehrMatignon2006}.

\medskip
{\bf Structure of the paper:} 
In Section \ref{SMC} we recall the definition and properties of semistable and stable marked curves which are relevant for us.
In Section \ref{SMC0} we review some special results in the case $g=0$ and give an explicit description of a Zariski neighborhood of a double point. We also prove a lemma about the extension of invertible sheaves that we use later on.
The stable marked model $(\CC,\CP_1,\ldots,\CP_{2g+2})$ of a hyperelliptic curve is constructed in Section \ref{HEC}, which is the core of this paper. 
In Section \ref{CF} we establish some properties of its special fiber and show how it can be computed directly from the special fiber of $(\bar\CC,\bar\CP_1,\ldots,\bar\CP_{2g+2})$. We also derive some consequences concerning the reduction of the Jacobian variety of~$C$.
In the final Section~\ref{Ex} we depict all combinatorial possibilities in the cases $g=1,2$ and exhibit further examples.

\section{Stable marked curves}
\label{SMC}

Throughout this paper, we consider a discrete valuation ring $R$ with quotient field $K$ and residue field $k=R/\Fm$ of characteristic $\neq 2$. We fix a uniformizer $t$ of $R$ and abbreviate $S := \Spec R$. 
In this section we review basic known facts about stable marked curves over~$S$. See Knudsen \cite{Knudsen1983} or Liu \cite[\S10.3]{LiuAlgGeo2002} for the general definition and properties of semistable and stable curves over arbitrary schemes.

Let $C$ be a geometrically connected smooth projective algebraic curve of genus $g$ over~$K$. A \emph{semistable model of~$C$} is a flat projective morphism $\CC\onto S$ with generic fiber~$C$, whose special fiber $C_0$ is a reduced curve that is smooth except possibly for finitely many ordinary double points. Each double point $p\in C_0$ possesses an \'etale neighborhood in $\CC$ which is \'etale over $\Spec R[x,y]/(xy-t^n)$ for some integer $n\ge1$, such that $p$ corresponds to the point $x=y=t=0$. Here the integer $n$ is uniquely determined by the local ring of $\CC$ at~$p$, for instance by Liu \cite[\S10.3]{LiuAlgGeo2002}. Following Liu \cite[\S5 Def.\,2]{Liu1993} we call $n$ the \emph{thickness of~$p$}.

Take distinct $K$-rational points $P_1,\ldots,P_n \in C(K)$. This turns $C$ into a \emph{smooth semistable marked curve} $(C,P_1,\ldots,P_n)$ over~$K$. If these points extend to pairwise disjoint sections $\CP_1,\ldots,\CP_n \in \CC(S)$ which avoid all double points of the special fiber, this turns $\CC$ into a \emph{semistable marked curve} $(\CC,\CP_1,\ldots,\CP_n)$ over~$S$. 

A \emph{stable marked curve} is a semistable marked curve such that for every fiber, the group of automorphisms which preserve the given sections is finite. In our situation this condition is satisfied for the smooth generic fiber if and only if $2g+n\ge3$. 
Granted this, the condition holds for the special fiber if and only if every non-singular rational geometric irreducible component contains at least $3$ marked or double points.
If $(\CC,\CP_1,\ldots,\CP_n)$ is stable, it is called a \emph{stable model of} $(C,P_1,\ldots,P_n)$, and its special fiber $(C_0,P_{0,1},\ldots,P_{0,n})$ is called a \emph{stable reduction of} $(C,P_1,\ldots,P_n)$.

Assume that the generic fiber $(C,P_1,\ldots,P_n)$ is stable. Then, after possibly replacing $R$ by a finite ramified extension, there exists a stable model, and it is unique up to unique isomorphism, once it exists (see for instance Liu \cite[2.19-21]{Liu2006} or Cuzub \cite[Th.\,3.4]{Cuzub2018}). 

\section{Stable marked curves of genus $0$}
\label{SMC0}

In this section we consider a geometrically connected smooth projective algebraic curve $\bar C$ of genus $0$ over $K$ that is marked by $n\ge3$ rational points $\bar P_1,\ldots,\bar P_n$. For any distinct indices $i$, $j$, $k\in\{1,\ldots,n\}$ there is then a unique isomorphism $(\bar C,\bar P_i,\bar P_j,\bar P_k) \cong (\BP^1_K,0,1,\infty)$. 
In this case it is known that, without extending~$R$, the stable marked curve $(\bar C,\bar P_1,\ldots,\bar P_n)$ possesses a stable model $(\bar \CC,\bar \CP_1,\ldots,\bar \CP_n)$ over~$S$, which is unique up to unique isomorphism. In fact, by the valuative criterion of properness this is a direct consequence of the main result of Gerritzen, Herrlich, and van der Put \cite[\S3]{GerritzenEtAl1988} that the functor of stable curves of genus $0$ with $n$ marked points is represented by a scheme which is smooth and projective over $\Spec\BZ$.

Moreover, this stable model can be computed explicitly by induction on~$n$. Namely, the stable model of $(\bar C,\bar P_1,\bar P_2,\bar P_3) \cong (\BP^1_K,0,1,\infty)$ is isomorphic to $(\BP^1_S,0,1,\infty)$, and for ${n>3}$ the stable model of $(\bar C,\bar P_1,\ldots,\bar P_n)$ is obtained from a stable model of $(\bar C,\bar P_1,\ldots,\bar P_{n-1})$ by the stabilization procedure of Knudsen \cite[Th.\,2.4]{Knudsen1983}. Each such step introduces at most one additional irreducible component isomorphic to $\BP^1_k$ in the special fiber. In particular this shows that every irreducible component of the special fiber $\bar C_0$ of $\bar\CC$ is isomorphic to~$\BP^1_k$. 

The procedure also shows that $\bar\CC$ can be covered by rational local charts, which makes computations easy. Alternatively such charts can be produced using the contraction process of \cite[Def.\,1.3, Prop.\,2.1]{Knudsen1983} or \cite[Lemma\,5]{GerritzenEtAl1988}. More explicitly:

\begin{Prop}\label{LocalChartsGenus0a}
Any smooth point $\bar p$ of $\bar C_0$ possesses a Zariski neighborhood in $\bar\CC$ which is isomorphic to an open subscheme of $\BA^1_S = \Spec R[x]$. 
If $\bar p$ is the reduction of a section~$\bar\CP_i$, this isomorphism can be chosen such that $\bar\CP_i$ corresponds to the section $x=0$.
\end{Prop}

\begin{Proof}
If $\bar C_0$ is irreducible, for any selection of distinct indices $i$, $j$, $k\in\{1,\ldots,n\}$ the tuple $(\bar\CC,\bar\CP_i,\bar\CP_j,\bar\CP_k)$ is a stable marked curve and therefore isomorphic to $(\BP^1_S,0,1,\infty)$,
proving everything. 

If $\bar C_0$ is reducible, the stability implies that $n>3$. Also, since the dual graph of $\bar C_0$ is a tree, there exists an irreducible component $\bar Z$ not containing $\bar p$ which meets the others in precisely one double point. By stability this component must meet some marked section~$\bar\CP_\ell$. Contracting $(\bar\CC,\bar\CP_1,\ldots,\bar\CP_{\ell-1},\bar\CP_{\ell+1},\ldots,\bar\CP_n)$ to a stable curve is then an isomorphism outside~$\bar Z$, and the statements follow by induction on~$n$.
\end{Proof}

\medskip

Next let $\CX$ be the flat projective curve over $S$ with sections $\CQ_1$, $\CQ_2$, whose special fiber consists of two irreducible components $Z_1$, $Z_2$, such that $\CX$ is covered by open charts 
$$\begin{array}{rll}
\CV_1 &\!=\ \CX \setminus(\CQ_2\cup Z_2) &\!=\ \Spec R[x^{-1}], \\[3pt]
\CV_2 &\!=\ \CX \setminus(\CQ_1\cup Z_1) &\!=\ \Spec R[y^{-1}], \\[3pt]
\CV_{12} &\!=\ \CX \setminus(\CQ_1\cup\CQ_2) &\!=\ \Spec R[x,y]/(xy-t^n)
\end{array}$$
for some integer $n>0$ and that
$$\begin{array}{l}
\hbox{$\CQ_1$ is the section of $\CV_1$ defined by $x^{-1}=0$,} \\[3pt]
\hbox{$\CQ_2$ is the section of $\CV_2$ defined by $y^{-1}=0$,} \\[3pt]
\hbox{$Z_1\cap\CV_{12}$ is defined by $y=t=0$, and} \\[3pt]
\hbox{$Z_2\cap\CV_{12}$ is defined by $x=t=0$.}
\end{array}$$
Note that the only singular point of the special fiber of~$\CX$ is the point $q$ of $\CV_{12}$ defined by $x=y=t=0$. The special fiber of $\CX$ thus looks as follows:
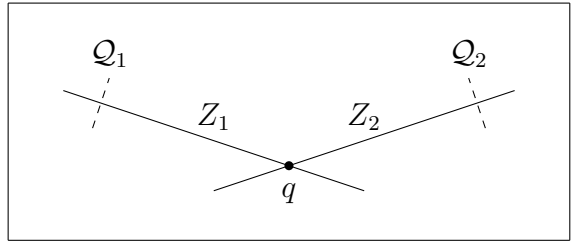
\[\FigZ\]

\begin{Prop}\label{LocalChartsGenus0}
For any double point $\bar p$ of $\bar C_0$ of thickness~$n$, there exists a morphism $\phi\colon\bar\CC\onto\CX$ which is an isomorphism near~$\bar p$, such that $\bar p$ maps to $q$ and that $\CQ_1$, $\CQ_2$ are the images of some marked sections $\bar\CP_i$, $\bar\CP_j$.
\end{Prop}

\begin{Proof}
Since the dual graph of $\bar C_0$ is a tree, the complement $\bar C_0\setminus\{\bar p\}$ consists of two connected components. By stability each of these must meet at least $2$ marked sections, say $\bar\CP_i$, $\bar\CP_{i'}$, respectively $\bar\CP_j$,~$\bar\CP_{j'}$. 
If $\bar p$ is the only double point of~$\bar C_0$, the tuple $(\bar\CC,\bar\CP_i,\bar\CP_{i'},\bar\CP_j,\bar\CP_{j'})$ is a stable marked curve with smooth generic fiber and singular special fiber. It well-known and implicit in \cite{GerritzenEtAl1988} that any such curve is isomorphic to $(\CX,\CQ_1,\CQ_3,\CQ_2,\CQ_4)$, where $(\CX,\CQ_1,\CQ_2)$ is as above except that $\CV_{12} = \Spec R[x,y]/(xy-s)$ for some $s\in\Fm\setminus\{0\}$, and where
 $$\begin{array}{l}
\hbox{$\CQ_3$ is the section of $\CV_1$ defined by $x^{-1}=1$, and} \\[3pt]
\hbox{$\CQ_4$ is the section of $\CV_2$ defined by $y^{-1}=1$.}
\end{array}$$
Writing $s=ut^n$ with $u\in R^\times$ and $n>0$ shows that $n$ is the thickness of~$\bar p$. After rescaling $y$ by $u$ we thus obtain the proposition with an isomorphism~$\phi$.

Now consider the case that there are other double points besides~$\bar p$. Since the dual graph of $\bar C_0$ is a tree, there exists an irreducible component $\bar Z$ not containing $\bar p$ which meets the others in precisely one double point. By stability this component must meet some marked section~$\bar\CP_\ell$. Contracting $(\bar\CC,\bar\CP_1,\ldots,\bar\CP_{\ell-1},\bar\CP_{\ell+1},\ldots,\bar\CP_n)$ to a stable curve is then an isomorphism outside~$\bar Z$, and the proposition follows by induction on~$n$.
\end{Proof}


\begin{Lem}\label{J*Lem}
Let $j$ denote the open embedding $\CV_{12}^\sm := \CV_{12}\setminus\{q\} \into \CV_{12}$ and consider an integer $0\le k<n$. 
\begin{enumerate}
\item[(a)] Then $j_*\CO_{\CV_{12}^\sm}(-kZ_2)$ is the coherent sheaf on $\CV_{12}$ associated to the ideal $(x,t^k)$ of $R[x,y]/(xy-t^n)$. 
\item[(b)] In particular $j_*\CO_{\CV_{12}^\sm} = \CO_{\CV_{12}}$.
\item[(c)] The module of relations between the generators of the ideal $(x,t^k)$ is generated by the relations $t^k\cdot x = x\cdot t^k$ and $y\cdot x = t^{n-k}\cdot t^k$.
\end{enumerate}
\end{Lem}

\begin{Proof}
By construction $\CV_{12}^\sm$ is covered by the open charts
$$\begin{array}{rl}
\CV_{12}\setminus Z_2 &\!=\ \Spec R[x^{\pm1}], \\[3pt]
\CV_{12}\setminus Z_1 &\!=\ \Spec R[y^{\pm1}],
\end{array}$$
on which the invertible sheaf $\CO_{\CV_{12}^\sm}(-kZ_2)$ is generated by~$1$, respectively by~$t^k$. Note that the rings $R[x^{\pm1}]$, $R[y^{\pm1}]$, $R[x,y]/(xy-t^n)$ all become subrings of $K[x^{\pm1}]$ on setting $y=\smash{\tfrac{t^n}{x}}$. Thus (a) amounts to the equality
$$R[x^{\pm1}] \;\cap\; R[(\tfrac{t^n}{x})^{\pm1}]\cdot t^k \ =\ R[x,\tfrac{t^n}{x}]\cdot x + R[x,\tfrac{t^n}{x}]\cdot t^k$$
within $K[x^{\pm1}]$. To prove this, we observe that both sides are graded for the $\BZ$-grading on $K[x^{\pm1}]$ which is trivial on~$K$ and for which $x$ is homogeneous of degree~$1$. For each of the subrings in question, the homogeneous part of degree $d\in\BZ$ is given by
$$\begin{array}{cl}
R[x^{\pm1}]_d &\!=\ R x^d, \\[6pt]
R[(\tfrac{t^n}{x})^{\pm1}]_d &\!=\ R x^d t^{-nd}, \\[3pt]
R[x,\tfrac{t^n}{x}]_d &\!=\ 
\biggl\{\begin{array}{ll}
R x^d &\hbox{if $d\ge0$,}\\[3pt]
R x^d t^{-nd} &\hbox{if $d\le0$.}
\end{array}
\end{array}$$
Thus the desired equality follows from the computations
$$\begin{array}{lcll}
R x^d \cap R x^d t^{-nd}\cdot t^k\ =\!& 
R x^d
\!&=\ R x^{d-1}\cdot x + R x^d\cdot t^k
&\hbox{if $d\ge1$,} \\[3pt]
R x^d \cap R x^d t^{-nd}\cdot t^k\ =\!& 
R x^d t^{-nd}\cdot t^k 
\!&=\ R x^{d-1} t^{-n(d-1)}\cdot x + 
R x^d t^{-nd} \cdot t^k
&\hbox{if $d\le0$.}
\end{array}$$
This proves (a), and (b) is the special case $k=0$ of (a). 

To prove (c) it suffices to describe the homogeneous relations of all degrees $d\in\BZ$. For $d\ge1$ these must have the form $ax^{d-1}\cdot x = bx^d\cdot t^k$ with $a,b\in R$. This equation requires that $a=bt^k$; hence all these relations are multiples of the relation $t^k\cdot x = x\cdot t^k$.
For $d\le0$ the relations must have the form $ay^{1-d}\cdot x = by^{-d}\cdot t^k$ with $a,b\in R$. This requires that $at^n = ayx =bt^k$ and hence $b=t^{n-k}a$; so all these relations are multiples of the relation $y\cdot x = t^{n-k}\cdot t^k$.
\end{Proof}

\section{Hyperelliptic curves}
\label{HEC}

Keeping the notation of the preceding sections, for the rest of this paper we assume that $\charact(k)\neq 2$. We fix a hyperelliptic curve $C$ of genus $g$ over~$K$. While usually this requires that $g\ge2$, everything below also holds for $g=1$.

\medskip
By assumption $C$ comes with a double covering $\pi\colon C \onto \bar C$ of a rational curve~$\bar C$. Since $K$ has characteristic $\not=2$, this covering is only tamely ramified. By the Hurwitz formula it is therefore ramified at exactly $2g+2$ geometric points, where $g$ is the genus of~$C$. After replacing $K$ by a finite extension, if necessary, we assume that these geometric points are all defined over~$K$. Let $P_1,\ldots,P_{2g+2} \in C(K)$ be these points, and let $\bar P_1,\ldots,\bar P_{2g+2} \in \bar C(K)$ denote their images under~$\pi$. Since $2g+2\ge4$, both $(C,P_1,\ldots,P_{2g+2})$ and $(\bar C,\bar P_1,\ldots,\bar P_{2g+2})$ are stable marked curves. In particular $\bar C\cong\BP^1_K$. 

\medskip
Let $\sigma$ denote the covering involution of $C$ over~$\bar C$. 

\begin{Lem}\label{Lgen}
\begin{enumerate}
\item[(a)] There is a unique decomposition
$$\pi_*\CO_C = \CO_{\bar C} \oplus \CL$$
with an invertible sheaf $\CL$ on~$\bar C$, such that $\sigma$ acts by 
$1$ on $\CO_{\bar C}$ and by $-1$ on~$\CL$.
\item[(b)]
The multiplication in $\pi_*\CO_C$ yields an isomorphism 
$$\CL^{\otimes2} \isoto \CO_{\bar C}(-\bar D)$$
for the divisor $\bar D := \bar P_1+\dots+\bar P_{2g+2}$ on~$\bar C$. 
\end{enumerate}
\end{Lem}

\begin{Proof}
Since $2$ is invertible on~$\bar C$, the direct image sheaf $\pi_*\CO_C$ is the direct sum of two sheaves on which $\sigma$ acts by $1$ and~$-1$, respectively. As $C/\langle\sigma\rangle\cong\bar C$, the first of these is the structure sheaf $\CO_{\bar C}$. But $\pi_*\CO_C$ is a locally free coherent sheaf of rank~$2$, because $\pi\colon C\onto\bar C$ is a double covering of smooth curves. Thus the second summand is locally free of rank~$1$, proving (a).

Since $\sigma$ acts by $(-1)^2=1$ on $\CL^{\otimes2}$, multiplication induces a homomorphism $\CL^{\otimes2} \to \CO_{\bar C}$. To determine its image consider a local chart $\Spec K[x]$ of~$\bar C$, over which $\CL$ corresponds to the module $K[x]\cdot y$. Then the resulting local chart of $C$ is $\Spec K[x,y]/(y^2-f(x))$ for some $f(x)\in K[x]$. The smoothness of $C$ implies that $f$ is separable and its zeros are precisely the branch locus of $\pi$ in our chart. Thus $\CL^{\otimes2}$ maps isomorphically to the ideal sheaf of $\bar D$, proving (b).\end{Proof}

\medskip

Next let $(\bar\CC,\bar\CP_1,\ldots,\bar\CP_{2g+2})$ be the stable model over $S:=\Spec R$ of the marked curve $(\bar C,\bar P_1,\ldots,\bar P_{2g+2})$. Let $\bar\CC^\sm$ denote the smooth locus of $\bar\CC$ over~$S$, that is, the complement of all double points in the special fiber. By semistability each section $\bar\CP_i$ already lands in~$\bar\CC^\sm$.


\begin{Lem}\label{Lsm}
After possibly replacing $R$ by any ramified extension of degree~$2$, there exists a unique extension of $\CL$ to an invertible sheaf $\CL^\sm$ on $\bar\CC^\sm$
such that the isomorphism in Lemma \ref{Lgen} (b) extends to an isomorphism 
$$(\CL^\sm)^{\otimes2} \isoto \CO_{\bar\CC^\sm}(-\bar\CD)$$
for the divisor $\bar\CD := \bar\CP_1+\dots+\bar \CP_{2g+2}$ on~$\bar\CC^\sm$.
\end{Lem}

\begin{Proof}
Choose a Weil divisor $\bar D'$ on $\bar C$ such that $\CL \cong \CO_{\bar C}(\bar D')$. Since $\bar\CC$ is regular, the closure $\bar\CD'$ of $\bar D'$ defines an invertible sheaf $\CL' \cong \CO_{\bar\CC^\sm}(\bar\CD')$ on $\bar\CC^\sm$ which extends~$\CL$. 
An arbitrary extension of $\CL$ to an invertible sheaf on $\bar\CC^\sm$ then has the form $\CL^\sm = \CL'(\sum n_{\bar X}\bar X)$ with an integer $n_{\bar X}$ for each irreducible component $\bar X$ of the special fiber of~$\bar\CC^\sm$. 

Now observe that the isomorphism in Lemma \ref{Lgen} (b) identifies $\CL^{\prime\otimes2}$ with an invertible sheaf on $\bar\CC^\sm$ which coincides with $\CO_{\bar C}(-\bar D)$ over~$\bar C$. Thus it induces an isomorphism
$$\CL^{\prime\otimes2}\ \isoto\ \CO_{\bar\CC^\sm}\bigl(-\bar\CD+{\textstyle\sum m_{\bar X}\bar X}\bigr)$$
with certain integers $m_{\bar X}$. For $\CL^\sm = \CL'(\sum n_{\bar X}\bar X)$ it therefore induces an isomorphism
$$(\CL^\sm)^{\otimes2}\ \isoto\ \CO_{\bar\CC^\sm}\bigl(-\bar\CD+{\textstyle\sum (m_{\bar X}+2n_{\bar X})\bar X}\bigr).$$
The problem thus has the unique solution $n_{\bar X} = -m_{\bar X}/2$, provided that all $m_{\bar X}$ are even. But this can be achieved, because replacing $R$ by any ramified extension of degree $2$ multiplies each $m_{\bar X}$ by~$2$.
\end{Proof}

\medskip
Now consider the inclusion $j\colon\bar\CC^\sm\into\bar\CC$. Since $\bar\CC$ is a normal scheme and the complement of $\bar\CC^\sm$ is a point of codimension~$2$, the sheaf $j_*\CL^\sm$ is coherent 
(see EGA4 \cite[Cor.\,5.11.4]{EGA4})
and the natural map $\CO_{\bar\CC}\to j_*\CO_{\bar\CC^\sm}$ is an isomorphism
(see \cite[D\'ef.\,5.9.9,
Th.\,5.10.5]{EGA4}).
%
%
These statements also follow from the explicit computations in the proof of Lemma \ref{CCdp} and from Lemma \ref{J*Lem} (b). The homomorphism in Lemma \ref{Lsm} therefore induces a homomorphism of sheaves of $\CO_{\bar\CC}$-modules
\UseTheoremCounterForNextEquation
\begin{equation}\label{jLsm}
(j_*\CL^\sm)^{\otimes2} \ \longto\ j_*\CO_{\bar\CC^\sm}(-\bar\CD)
\ =\ \CO_{\bar\CC}(-\bar\CD)\ \longinto\ \CO_{\bar\CC}.
\end{equation}
This turns $\CO_{\bar\CC} \oplus j_*\CL^\sm$ into a coherent sheaf of $\CO_{\bar\CC}$-algebras. 
We are interested in the relative $\Spec$ and its canonical finite morphism
\UseTheoremCounterForNextEquation
\begin{equation}\label{CCPiDef}
\Pi\colon\ \CC := \Spec(\CO_{\bar\CC} \oplus j_*\CL^\sm)\ \longonto\ \bar\CC.
\end{equation}
By construction $\CC$ has generic fiber~$C$, and $\Pi$ extends the given double covering $\pi\colon C\onto \bar C$. 


\begin{Lem}\label{CCsm}
\begin{enumerate}
\item[(a)] The inverse image $\CC^\sm := \Pi^{-1}(\bar\CC^\sm)$ is a smooth curve over~$S$.
\item[(b)] Each section $\bar\CP_i$ lifts to a unique section $\CP_i\in\CC^\sm(S)$ which extends $P_i\in C(K)$. 
\item[(c)] The morphism $\CC^\sm \to \bar\CC^\sm$ is finite of degree~$2$ and ramified exactly along the~$\CP_i$.
\end{enumerate}
\end{Lem}

\begin{Proof}
The first statement in (c) follows from the fact that $\CO_{\bar\CC^\sm} \oplus \CL^\sm$ is locally free of rank~$2$.
To prove the remaining statements, by Proposition \ref{LocalChartsGenus0a} we can restrict ourselves to an open chart $\CU \subset \bar\CC^\sm$ that is isomorphic to an open subscheme of $\BA^1_S = \Spec R[x]$. Then the restriction $\CL^\sm|\kern2pt\CU$ is free, say with generator~$u$, and $\Pi^{-1}(\CU)$ is an open subscheme of $\Spec R[x,u]/(u^2-f)$ for some $f\in R[x]$. By Lemma \ref{Lsm} the zero locus of $f$ on $\CU$ is precisely the divisor $\CU\cap\bar\CD$. 

Over $\CU\setminus\bar\CD$ the equation $u^2=f$ is therefore \'etale and $\Pi^{-1}(\CU\setminus\bar\CD)$ is a smooth curve. Near a point of $\bar\CP_i$, by Proposition \ref{LocalChartsGenus0a} we can assume that $\bar\CP_i$ is given by the equation $x=0$. Then $f=xg$ for some $g\in R[x]$ with a unit $g(0)\in R^\times$. Thus $\frac{df}{dx}(0) = g(0)\in R^\times$; hence the equation $u^2=f$ defines a smooth curve near $\bar\CP_i$ with the local parameter~$u$, and the morphism $\Pi$ is ramified over~$\bar\CP_i$. Finally, the unique lift of $\bar\CP_i$ is given by the local equations $x=u=0$. 
\end{Proof}

\medskip

Now let $\bar p$ be a double point of~$\bar C_0$. Since the dual graph of $\bar C_0$ is a tree, the complement $\bar C_0\setminus\{\bar p\}$ consists of two connected components. This divides the $2g+2$ sections $\bar\CP_i$ into two groups.

\begin{Def}\label{EvenOddDef}
We call $\bar p$ \emph{even} (resp.\ \emph{odd}), if each connected component of $\bar C_0\setminus\{\bar p\}$ meets an even (resp.\ odd) number of sections~$\bar\CP_i$.
\end{Def}

\begin{Lem}\label{CCdp}
\begin{enumerate}
\item[(a)] The scheme $\CC$ is a semistable curve over~$S$.
\item[(b)] If $\bar p$ is even, the morphism $\Pi\colon\CC\onto\bar\CC$ is \'etale over~$\bar p$, and the inverse image of $\bar p$ has two geometric points which are double points of the same thickness as~$\bar p$. 
\item[(c)] If $\bar p$ is odd, the inverse image of $\bar p$ is a double point whose thickness is half the thickness of~$\bar p$.
\end{enumerate}
\end{Lem}

\begin{Proof}
Since $\CC^\sm := \Pi^{-1}(\bar\CC^\sm)$ is already a smooth curve by Lemma \ref{CCsm} (a), it suffices to prove everything near a double point $\bar p$ of~$\bar C_0$. For this choose a morphism $\phi\colon\bar\CC\onto\CX$ with $\phi( \bar p)= q$ as in Proposition \ref{LocalChartsGenus0}. Then the special fibers of $\phi^{-1}(\CV_1)$ and $\phi^{-1}(\CV_2)$ are the two connected components of $\bar C_0\setminus\{\bar p\}$. For each $i$ abbreviate $\CR_i := \phi(\bar\CP_i)$. After renumbering the sections we may without loss of generality assume that $\CR_1,\ldots,\CR_r$ land in $\CV_1$ and $\CR_{r+1},\ldots,\CR_{2g+2}$ land in~$\CV_2$,
and that $\CR_1=\CQ_1$ and $\CR_{2g+2}=\CQ_2$. Then $\bar p$ is even (resp.\ odd) if and only if $r$ is even (resp.\ odd).

Suppose that $\CR_i$ is given by the equation $x^{-1}=\lambda_i \in R$ on $\CV_1=\Spec R[x^{-1}]$ for $i\le r$, respectively by the equation $y^{-1}=\mu_i \in R$ on $\CV_2=\Spec R[y^{-1}]$ for $i>r$. Then 
$$f\ :=\ \prod_{i=1}^r (1-\lambda_ix) \cdot \!\!\prod_{i=r+1}^{2g+2} (1-\mu_iy)$$
is an element of $R[x,y]/(xy-t^n)$ with zero divisor $\phi(\bar\CD)\cap\CV_{12}$ that is invertible at~$q$. 

Now observe that $\phi\colon\bar\CC\onto\CX$ is an isomorphism outside a finite closed subset of $\CX$ not containing~$q$. Since invertible sheaves on regular schemes correspond to divisor classes, the sheaf $\CL^\sm$ corresponds to an invertible sheaf on $\CX\setminus\{q\}$.
For simplicity we denote this sheaf again by $\CL^\sm$.
Over $\CV_1 = \Spec R[x^{-1}]$ it is free, say with generator~$u_1$. The homomorphism from Lemma \ref{Lsm} then sends $u_1^2$ to an element of $R[x^{-1}]$ with zero divisor $\phi(\bar\CD)\cap\CV_1$. Our description of the points $\CR_i$ and the relation $xy=t^n$ thus implies that
\UseTheoremCounterForNextEquation
\begin{equation}\label{u2Rel}
u_1^2\ =\ c\cdot \prod_{i=1}^r (x^{-1}-\lambda_i) \cdot\!\! \prod_{i=r+1}^{2g+2} (1-\mu_i t^nx^{-1})
\ =\ c\cdot x^{-r}\cdot f
\end{equation}
for some constant $c\in R^\times$. Likewise, the sheaf $\CL^\sm$ is free on $\CV_2 = \Spec R[y^{-1}]$, say with generator~$u_2$, and
\UseTheoremCounterForNextEquation
\begin{equation}\label{v2Rel}
u_2^2\ =\ d \cdot \prod_{i=1}^r (1-\lambda_i t^n y^{-1}) \cdot \!\!\prod_{i=r+1}^{2g+2} (y^{-1}-\mu_i)
\ =\ d\cdot y^{r-2g-2}\cdot f
\end{equation}
for some constant $d\in R^\times$. Together this implies that
$$u_1^2\ =\ \frac{c}{x^r}\cdot f
\ =\ \frac{c}{x^r}\cdot \frac{y^{2g+2-r}}{d}\cdot u_2^2
\ =\ \frac{c}{d}\cdot \frac{y^{2g+2}}{t^{nr}}\cdot u_2^2$$
Therefore $nr$ must be even, and we have
\UseTheoremCounterForNextEquation
\begin{equation}\label{uvRel}
u_1\ =\ e\cdot\frac{y^{g+1}}{t^{nr/2}}\cdot u_2
\end{equation}
for some constant $e\in R^\times$ with $e^2=\frac{c}{d}$.

\medskip
We can now compute $j_*\CL^\sm$ near the double point~$q$. Recall that the deleted neighborhood $\CV_{12}^\sm = \CV_{12}\setminus\{q\}$ of $q$ is covered by the affine open charts $\CV_1\cap\CV_{12} = \Spec R[x^{\pm1}]$ and $\CV_2\cap\CV_{12} = \Spec R[y^{\pm1}]$. Since $x$ is invertible on $\CV_1\cap\CV_{12}$, the restriction of $\CL^\sm$ to $\CV_1\cap\CV_{12}$ is generated by any element of the set $R^\times x^\BZ u_1$. Likewise, the restriction of $\CL^\sm$ to $\CV_2\cap\CV_{12}$ is generated by any element of $R^\times y^\BZ u_2$.

\medskip
Suppose first that $r$ is even. Using $xy=t^n$ the equation (\ref{uvRel}) then implies that
$$u\ :=\ x^{r/2}\cdot u_1\ =\ e\cdot y^{g+1-r/2}\cdot u_2.$$
By construction this is a generator of $\CL^\sm$ over both $\CV_1\cap\CV_{12}$ and $\CV_2\cap\CV_{12}$. Thus $\CL^\sm$ is free with basis $u$ in a deleted neighborhood of~$q$. The equality $j_*\CO_{\CV_{12}^\sm} = \CO_{\CV_{12}}$ from Lemma \ref{J*Lem} (b) now implies that $j_*\CL^\sm$ is free with basis $u$ in a neighborhood of~$q$. Moreover, since the divisor $\bar\CD$ does not meet~$\bar p$, the homomorphism
$$(j_*\CL^\sm)^{\otimes2} \ \longto\ j_*\CO_{\CV_{12}^\sm} \ =\ \CO_{\CV_{12}}$$
corresponding to (\ref{jLsm}) is an isomorphism in a neighborhood of~$q$, where the last equality again follows from Lemma \ref{J*Lem} (b). In fact, we already know that $u^2 = x^r u_1^2 = cf$ is a unit near~$q$. 
Thus $\CC$ is locally over $\bar p$ given by the \'etale equation $u^2=cf$. Since $cf$ is invertible at~$q$, this proves (b).

\medskip
Suppose now that $r$ is odd. Then $n$ must be even, and using $xy=t^n$ the equation (\ref{uvRel}) implies that
$$u\ :=\ x^{(r-1)/2}\cdot u_1\ =\ \frac{e\cdot y^{g+1-(r-1)/2}\cdot u_2}{t^{n/2}}.$$
By construction $\CL^\sm$ is generated by $u$ over $\CV_1\cap\CV_{12}$ and by $t^{n/2}u$ over $\CV_2\cap\CV_{12}$. Thus multiplication by $u$ induces an isomorphism 
$\CO_{\CV_{12}^\sm}(-\frac{n}{2} Z_2) \stackrel{\sim}{\longto} \CL^\sm|\kern2pt\CV_{12}^\sm$, and therefore also an isomorphism $j_*\CO_{\CV_{12}^\sm}(-\frac{n}{2}Z_2) \stackrel{\sim}{\longto} j_*\CL^\sm|\kern2pt\CV_{12}$. By Lemma \ref{J*Lem} (a) we conclude that $j_*\CL^\sm$ is the coherent sheaf on $\CV_{12}$ that is associated to the $R[x,y]/(xy-t^n)$-module $(x,t^{n/2})\cdot u$. 
Consider the generators $v := xu$ and $w := t^{n/2}u$ of this module. Lemma \ref{J*Lem} (c) implies that the module of relations between them is generated by the relations 
$$\begin{array}{rl}
t^{n/2} v &\!=\ x w, \\[3pt]
y v  &\!=\ t^{n/2} w.
\end{array}$$
Using (\ref{u2Rel}) we also obtain the relations
$$\begin{array}{rl}
v^2 &\!=\ x^{r+1} u_1^2     \ =\ xcf, \\[3pt]
w^2 &\!=\ t^nx^{r-1} u_1^2  \ =\ ycf, \\[3pt]
vw  &\!=\ t^{n/2} x^r u_1^2 \ =\ t^{n/2}cf.
\end{array}$$
Together this shows that over a neighborhood of $\bar p$ the scheme $\CC$ is open in the spectrum of the ring
\UseTheoremCounterForNextEquation
\begin{equation}\label{APresent}
A\ :=\ R\bigl[x,y,v,w\bigr]\bigm/
\Bigl(\begin{array}{c}
xy-t^n,\ t^{n/2}v-xw,\ yv-t^{n/2}w, \\
v^2-xcf,\ w^2-ycf,\ vw-t^{n/2}cf
\end{array}\Bigr).
\end{equation}
Recall that $\bar p$ corresponds to the point $x=y=t=0$. Thus the equations $v^2=xcf$ and $w^2=ycf$ show that $\CC$ has a unique point $p$ defined by ${x=y=t=v=w=0}$ above~$\bar p$. Moreover, since $cf$ is invertible at~$p$, the first three relations defining $A$ are linear combinations of the last three in a neighborhood of~$p$. Furthermore, in the completion of $A$ at $p$ we can eliminate the variables $x$ and $y$ using the equations $v^2=xcf$ and $w^2=ycf$. The completion is therefore isomorphic to 
$$R[[v,w]]/(vw-t^{n/2}cf)$$
with $cf\in R[[v,w]]^\times$. This shows that $\CC$ is a semistable curve near $p$ and that $p$ is a double point whose thickness $n/2$ is half the thickness $n$ of~$\bar p$, proving (c).
\end{Proof}


\begin{Prop}\label{CCstable}
$(\CC,\CP_1,\ldots,\CP_{2g+2})$ is a stable model of $(C,P_1,\ldots,P_{2g+2})$.
\end{Prop}

\begin{Proof}
By Lemmas \ref{CCsm} and \ref{CCdp} the scheme $\CC$ is a semistable model over~$S$ and the sections $\CP_i$ land in the smooth part~$\CC^\sm$. Moreover, since the $\CP_i$ lift the disjoint sections $\bar\CP_i$ of~$\bar\CC$, they are themselves disjoint. Therefore $(\CC,\CP_1,\ldots,\CP_{2g+2})$ is a semistable model of $(C,P_1,\ldots,P_{2g+2})$.

To prove the stability we first observe that the construction is invariant under arbitrary finite extensions of~$R$. Thus we may assume that every irreducible component $Z$ of the special fiber $C_0$ of~$\CC$ is geometrically irreducible. For any such $Z$ the image $\Pi(Z)$ is an irreducible component of~$\bar C_0$. By the stability of $(\bar\CC,\bar\CP_1,\ldots,\bar\CP_{2g+2})$ this irreducible component contains at least $3$ marked or double points. But by Lemmas \ref{CCsm} and \ref{CCdp}, each point of $C_0$ above a marked or double point of $\bar C_0$ is itself a marked or double point. Therefore $\bar Z$ contains at least $3$ marked or double points, and $(\CC,\CP_1,\ldots,\CP_{2g+2})$ is stable.
\end{Proof}

\begin{Rem}\label{CCQuot}
\rm The hyperelliptic involution $\sigma$ on $C$ extends uniquely to an involution of $\CO_{\bar\CC} \oplus j_*\CL^\sm$ and hence of~$\CC$ by letting $\sigma$ act on $\CO_{\bar\CC}$ by $1$ and on $j_*\CL^\sm$ by $-1$. This construction directly implies that $\Pi$ induces an isomorphism $\CC /\langle\sigma\rangle \cong \bar\CC$. 
By the uniqueness of stable models it follows that for any stable model of $(C,P_1,\ldots,P_{2g+2})$, the quotient by $\langle\sigma\rangle$ is a stable model of $(\bar C,\bar P_1,\ldots,\bar P_{2g+2})$. In this form, however, we do not have a simple direct proof for the statement.
\end{Rem}

\section{The special fiber}
\label{CF}

The following diagram collects the schemes and sections constructed in Section~\ref{HEC}. In addition $\pi_0\colon C_0\onto\bar C_0$ denotes the morphism of special fibers induced by~$\Pi$,
and $P_{0,i} \in C_0(k)$ and $\bar P_{0,i}\in\bar C_0(k)$ denote the respective induced rational points therein:
$$\xymatrix{
\ C\  \ar@{^{ (}->}[r] \ar@{->>}[d]_-\pi & 
\ \CC\ \ar@{->>}[d]_-\Pi & 
\ C_0\ \ar@{_{ (}->}[l] \ar@{->>}[d]_-{\pi_0} &&
\ P_i\ \ar@{^{ (}->}[r] \ar@{|->}[d] & 
\ \CP_i\ \ar@{|->}[d] & 
\ P_{0,i}\ \ar@{|->}[d] \ar@{_{ (}->}[l] \\
\ \bar C\ \ar@{^{ (}->}[r] \ar@{->>}[d] & 
\ \bar\CC\ \ar@{->>}[d] & 
\ \bar C_0\ \ar@{->>}[d] \ar@{_{ (}->}[l] &&
\ \bar P_i\ \ar@{^{ (}->}[r] & 
\ \bar\CP_i\ & 
\ \bar P_{0,i}\ \ar@{_{ (}->}[l] \\
\ \Spec K\ \ar@{^{ (}->}[r] & 
\ S\  & 
\ \Spec k\ \ar@{_{ (}->}[l] \\}$$
In this section we analyze more closely the special fiber $(C_0,P_{0,1},\ldots,P_{0,n})$.


\begin{Prop}\label{IrredComp}
Let $\bar Z$ be an irreducible component of~$\bar C_0$, and let $Z := \pi_0^{-1}(\bar Z)^\red$ be the reduced closed subscheme of its inverse image. Let $T$ be the set of marked points $\bar P_{0,i}$ and odd (!) double points of~$\bar C_0$ which lie in~$\bar Z$.
\begin{enumerate}
\item[(a)] Then $Z$ is smooth, and the morphism $Z \onto \bar Z$ is a double covering that is ramified precisely over the points in~$T$.
\item[(b)] If $T=\emptyset$, then $Z$ is isomorphic to the disjoint union of two copies of~$\BP^1_k$ after possibly replacing $R$ by an unramified extension of degree~$2$.
\item[(c)] If $T\not=\emptyset$, then $Z$ is a geometrically irreducible component of~$C_0$, and its genus $g_Z$ satisfies $|T|=2g_Z+2$.
In particular $|T|$ is even.
\end{enumerate}
\end{Prop}

\begin{Proof}
Since $\bar C_0$ is a semistable curve of genus~$0$, all its irreducible components are smooth. Thus (a) follows from Lemma \ref{CCsm} outside all double points of~$\bar C_0$, and from Lemma \ref{CCdp} (b) over any even double point. 
For any odd double point $\bar p\in T$ we choose a morphism $\phi\colon\bar\CC\onto\CX$ as in the proof of Lemma \ref{CCdp}. Without loss of generality we may assume that $\phi(\bar Z) = Z_1$, which is defined by the equations $y=t=0$ in the coordinate ring $R[x,y]/(xy-t^n)$ of~$\CV_{12}$. Dividing the ring (\ref{APresent}) by the ideal $(y,t)$ then shows that $\pi_0^{-1}(\bar Z)$ is locally over $\bar p$ isomorphic to the spectrum of the ring
$$k\bigl[x,v,w\bigr]\bigm/
\bigl(xw,\ v^2-x\cdot\overline{cf},\ w^2,\ vw\bigr),$$
where $\overline{cf}$ denotes the residue class of~$cf$. Here the residue class of $w$ is nilpotent, so dividing further by $(w)$ shows that $Z$ is locally over $\bar p$ isomorphic to the spectrum of the ring
$$k[x,v]/\bigl(v^2-x\cdot\overline{cf}\kern1pt\bigr)$$
modulo its nilradical. But since $cf$ is invertible at~$q$, the equation $v^2=x\cdot\overline{cf}$ already defines a smooth curve at $x=0$, which is ramified of degree $2$ over $Z_1\cap\CV_{12} = \Spec k[x]$. This proves (a) over any odd double point.

If $T=\emptyset$, then $Z\onto\bar Z\cong\BP^1_k$ is an everywhere unramified double covering by (a). Thus after possibly replacing $k$ by an extension of degree~$2$, it becomes isomorphic to the disjoint union of two copies of~$\BP^1_k$, proving (b).

If $T\not=\emptyset$, then $Z\onto\bar Z\cong\BP^1_k$ is ramified somewhere; hence $Z$ is geometrically irreducible. As $\charact(k)\not=2$, by the Hurwitz formula this covering is ramified at exactly $2g_Z+2$ geometric points. Since all points in $T$ are already defined over~$k$, it follows that $|T|=2g_Z+2$, proving~(c).
\end{Proof}


\begin{Prop}\label{DualGraph}
After possibly replacing $k$ by a finite extension, we have:
\begin{enumerate}
\item[(a)] The special fiber $(C_0,P_{0,1},\ldots,P_{0,2g+2})$ up to isomorphism can be computed directly from the special fiber $(\bar C_0,\bar P_{0,1},\ldots,\bar P_{0,2g+2})$.
\item[(b)] The dual graph of $C_0$ and its marking by the points $P_{0,1},\ldots,P_{0,2g+2}$ as well as the genus of each irreducible component can be computed directly from the dual graph of $\bar C_0$ and its marking by the points $\bar P_{0,1},\ldots,\bar P_{0,2g+2}$.
\end{enumerate}
\end{Prop}


\begin{Proof}
After replacing $k$ by a finite extension, for every irreducible component $\bar Z$ of $\bar C_0$ the covering $Z\onto\bar Z$ from Proposition \ref{IrredComp} (a) is determined up to isomorphism by the branch points. It is therefore determined by $(\bar C_0,\bar P_{0,1},\ldots,\bar P_{0,2g+2})$. This also determines the points $P_{0,1},\ldots,P_{0,2g+2}$. To prove (a) it remains to show that the gluing at the double points is determined up to isomorphism. 

For this recall that the dual graph of $\bar C_0$ is a tree. We can therefore order the irreducible components of $\bar C_0$ in such a way that each next one meets the union of the earlier ones at precisely one double point. Let $\bar Z_1,\ldots,\bar Z_m$ be these irreducible components and put $Z_i := \pi_0^{-1}(\bar Z_i)^\red$. 
For every $1\le i\le m$ set $\bar Z_{\le i} := \bar Z_1\cup\ldots\cup\bar Z_i$ and $Z_{\le i} := \pi_0^{-1}(\bar Z_{\le i})^\red$. 

Suppose that we have already constructed $Z_{\le i}$ for some $1\le i<m$. Then by assumption $\bar Z_{i+1}$ meets $\bar Z_{\le i}$ at precisely one double point~$\bar p$. If $\bar p$ is odd, by Proposition \ref{IrredComp} (a) both $Z_{i+1}$ and $Z_{\le i}$ have a unique rational point above~$\bar p$, and by Lemma \ref{CCdp} (c) we must glue them there, for which there is precisely one way.
If $\bar p$ is even, by Proposition \ref{IrredComp} (a) both $Z_{i+1}$ and $Z_{\le i}$ have two geometric points above~$\bar p$. After possibly replacing $k$ by a finite extension we may assume them to be rational. By Lemma \ref{CCdp} (b) we must glue them in pairs. But the two ways of pairing them are interchanged by the covering involution of $Z_{i+1}\onto\bar Z_{i+1}$. Thus up to isomorphism this pairing is unique, and then the gluing is unique. In both cases we can therefore construct $Z_{\le i}$ uniquely up to isomorphism.

By induction on $i$ this shows that we can construct $Z_{\le m}=Z$ uniquely up to isomorphism, proving (a).

\medskip
The same procedure, applied to only the combinatorial data in the dual graph of $\bar C_0$ and its marking by the points $\bar P_{0,1},\ldots,\bar P_{0,2g+2}$ without their precise location, yields the dual graph of $C_0$ and its marking by the points $P_{0,1},\ldots,P_{0,2g+2}$. Finally, this data also determines the genus of each irreducible component of $C_0$ by Proposition \ref{IrredComp}. This proves (b).
\end{Proof}


\begin{Lem}\label{Conn}
Let $\bar p$ be a double point of~$\bar C_0$, and let $\bar X$ be the closure of a connected component of $\bar C_0\setminus\{\bar p\}$. Then the inverse image $\pi_0^{-1}(\bar X)$ is connected.
\end{Lem}

\begin{Proof}
By stability $\bar X$ must contain a marked point~$\bar P_{0,j}$. Let $\bar Z_1$ be the irreducible component of $\bar C_0$ that contains~$\bar P_{0,j}$. In the proof of Proposition \ref{DualGraph} we can choose $\bar Z_2,\ldots,\bar Z_m$ by first exhausting all irreducible components of $\bar X$ and then passing through the double point $\bar p$ to collect the remaining irreducible components. Then $\bar X = \bar Z_{\le k}$ for some $1\le k<m$.

Now $Z_1$ is irreducible by Proposition \ref{IrredComp} (c), because $\bar Z_1$ contains a marked point. In particular $Z_1$ is connected. By induction on $i$ the gluing procedure in the proof of Proposition \ref{DualGraph} shows that $Z_{\le i}$ is connected for every~$i$. In particular $\pi_0^{-1}(\bar X) = Z_{\le k}$ is connected, as desired.
\end{Proof}

\medskip

Now we look at the Jacobian variety $J$ of~$C$. Recall from Bosch-L\"{u}tkebohmert-Raynaud \cite[\S9.7,\,Cor.\,2]{BoschRaynaudLuetkerbohmert1990} that whenever $C$ has semistable reduction~$C_0$, then so does $J$ and the identity component $J_0$ of its reduction is the identity component of the Picard scheme of~$C_0$. There is an integer $r$ such that $J_0$ is an  extension of an abelian variety of dimension $g-r$ by a torus of dimension~$r$.

Let $n$ denote the number of double points of $\bar C_0$ and $n_0$ the number of even double points. Let $m$ denote the number of irreducible components of $\bar C_0$ and $m_0$ the number of irreducible components which contain no marked point or odd double point.

\begin{Prop}\label{JacRed}
\begin{enumerate}
\item[(a)] We have $r=n_0-m_0$.
\item[(b)] We have $r=0$ if and only if $n_0=0$.
In other words $J$ has potential good reduction if and only if $\bar C_0$ possesses no even double point.
\end{enumerate}
\end{Prop}

\begin{Proof}
(For another proof see
\cite[Thm.\,10.3 (6)-(7)]{DDMM}.)
After extending $k$ we may without loss of generality assume that all double points of $C_0$ are defined over~$k$ and all irreducible components of $C_0$ are geometrically irreducible. Then by Lemma \ref{CCdp} every odd double point of $\bar C_0$ gives rise to one double point of~$C_0$, and every even double point of $\bar C_0$ gives rise to two double points. The number of double points of $C_0$ is therefore $n+n_0$. Likewise, by Proposition \ref{IrredComp} every irreducible component of $\bar C_0$ that contains a marked or odd double point gives rise to one irreducible component of~$C_0$, and the remaining irreducible components of $\bar C_0$ give rise to two irreducible components each. The number of irreducible components of $C_0$ is therefore $m+m_0$. 

The dual graph of $C_0$ thus has $m+m_0$ vertices and $n+n_0$ edges. Likewise the dual graph of $\bar C_0$ has $m$ vertices and $n$ edges. As this graph is a tree, we have $m=n+1$. As $C_0$ is connected, so is its dual graph. The $h^1$ 
of this graph is therefore ${(n+n_0)-(m+m_0)+1} = n_0-m_0$. By \cite[\S9.2,\,Ex.\,8]{BoschRaynaudLuetkerbohmert1990} this proves (a).

To prove (b) note first that if $n_0=0$, then (a) implies that $r\le0$ and hence $r=0$. So assume that $n_0>0$ and let $\bar p$ be an even double point of~$\bar C_0$. As the dual graph of $\bar C_0$ is a tree, the complement $\bar C_0\setminus\{\bar p\}$ consists of two connected components. Let $\bar X$ and $\bar X'$ denote their closures in~$\bar C_0$. Then both $\pi_0^{-1}(\bar X)$ and $\pi_0^{-1}(\bar X')$ are connected by Lemma \ref{Conn}. But in the proof of Proposition \ref{DualGraph} we have seen that they are glued together at two double points. Thus the dual graph of $C_0$ has a cycle and hence $r>0$, as desired.
\end{Proof}


\section{Examples}
\label{Ex}

Each of the following pictures shows the closed fiber $C_0$ above the closed fiber~$\bar C_0$, such that the morphism $\pi_0\colon C_0\to\bar C_0$ is approximately the vertical projection. Each solid line or curve represents an irreducible component and each intersection a double point. For simplicity we assume that all these are geometrically irreducible, resp.\ defined over~$k$. The even double points of $\bar C_0$ are enclosed in an additional small circle, and the pictures show that there are two double points of $C_0$ above each of them.
An irreducible component of genus $g'>0$ is labeled by~$g'$, while all irreducible components of genus~$0$ remain unlabeled. The dashed lines represent the $2g+2$ marked sections. 

\medskip

\begin{figure}[h]
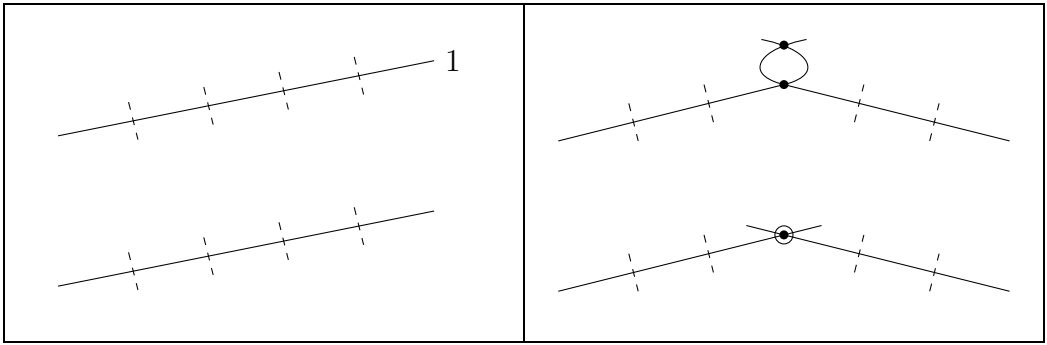

\centering
\FigA
\caption{The two possibilities in genus $g=1$}
\label{FigA}
\end{figure}

In the case $g=1$ there are exactly $2$ combinatorial possibilities for $(\bar C_0,\bar P_{0,1},\ldots,\bar P_{0,2g+2})$, which result in the description of $(C_0,P_{0,1},\ldots,P_{0,2g+2})$ in Figure~\ref{FigA}.

\medskip

\begin{figure}[p]
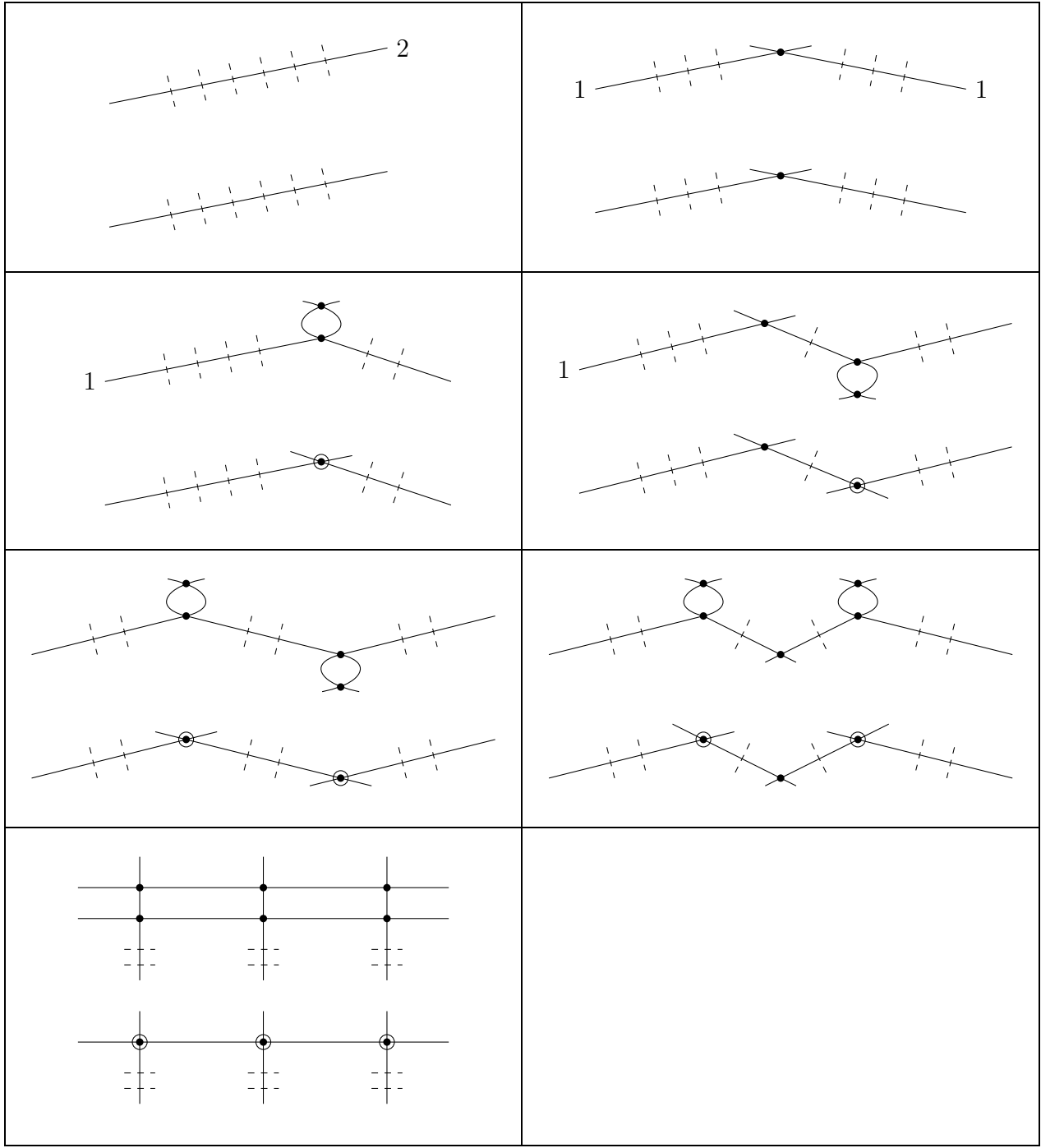

\centering
\FigB
\caption{The seven possibilities in genus $g=2$}
\label{FigB}
\end{figure}

In the case $g=2$ there are exactly $7$ cases shown in Figure~\ref{FigB}. In the last of these $\bar C_0$ has an irreducible component over which there are two irreducible components of~$C_0$.

\medskip

The example for $g=3$ in Figure~\ref{FigC} shows how the irreducible components of $C_0$ can intersect when there are many even double points.

\begin{figure}[h!]
\centering
\FigC
\caption{An example of genus $g=3$}
\label{FigC}
\end{figure}



\medskip

The example for $g=7$ in Figure~\ref{FigD} illustrates how different features from the other cases can occur simultaneously.

\begin{figure}[h!]
\centering
\FigD
\caption{An example of genus $g=7$}
\label{FigD}
\end{figure}



\end{document}

%% file: GP-Figures.tex
\def\FigZ{
\begin{tikzpicture}[
       outer frame sep=0ex,%
      /tikz/inner frame xsep=4ex,
      /tikz/inner frame ysep=2ex,
      show background top,%
      show background bottom,%
      show background left,%
      show background right]
      \draw plot (0,0) -- (4,1+1/3);
      \draw (7/2-1/9, 1+1/3-1/6+1/3) node [above, scale=1] {$\CQ_2$};
      \draw[dashed] plot  (7/2+1/9, 1-1/3-1/6+1/3) -- (7/2-1/9, 1+1/3-1/6+1/3);
      \draw (4/2, 1-1/3) node [above, scale=1] {$Z_2$};
      \begin{scope}[xscale=-1, shift={(-2,0)}]
       \draw plot (0,0) -- (4,1+1/3);
       \draw[dashed] plot  (7/2+1/9, 1-1/3-1/6+1/3) -- (7/2-1/9, 1+1/3-1/6+1/3);
       \draw (7/2-1/9, 1+1/3-1/6+1/3) node [above, scale=1] {$\CQ_1$};
       \draw (4/2, 1-1/3) node [above, scale=1] {$Z_1$};
      \end{scope}
\fill [black]   (1,1/3) circle (0.06);%
\draw (1,1/3-1/12) node [below, scale=1] {$q$};
\end{tikzpicture}
}

\def\FigA{
\begin{tabular}{|m{0.4\textwidth}|m{0.4\textwidth}|}
       \hline 
       \begin{center}
       \begin{tikzpicture}
              \draw plot (0,0) -- (5,1);
              \draw[dashed] plot  (2-1/16,2/5+1/4) -- (2+1/16,2/5-1/4);
              \draw[dashed] plot  (1-1/16,1/5+1/4) -- (1+1/16,1/5-1/4);
              \draw[dashed] plot  (3-1/16,3/5+1/4) -- (3+1/16,3/5-1/4);
              \draw[dashed] plot  (4-1/16,4/5+1/4) -- (4+1/16,4/5-1/4);
              \draw (5,1) node [right, scale=1] {1};
              \begin{scope}[shift={(0,-2)}]
                     \draw plot (0,0) -- (5,1);
                     \draw[dashed] plot  (2-1/16,2/5+1/4) -- (2+1/16,2/5-1/4);
              \draw[dashed] plot  (1-1/16,1/5+1/4) -- (1+1/16,1/5-1/4);
              \draw[dashed] plot  (3-1/16,3/5+1/4) -- (3+1/16,3/5-1/4);
              \draw[dashed] plot  (4-1/16,4/5+1/4) -- (4+1/16,4/5-1/4);
              \end{scope}
             
              \end{tikzpicture}
       \end{center}
              &
              \begin{center}
              \begin{tikzpicture}

             \draw plot (0,0) -- (3,6/8);
             \begin{scope}[shift={(3,3/4)}, scale=0.6]
              \draw plot[smooth, tension=1.3] coordinates {(0, 0) (0.5, 1/2) (-0.5, 1)};
              \draw plot[smooth, tension=1.3] coordinates {(0, 0) (-0.5, 1/2) (0.5, 1)};
              \fill [black]   (0,7/8) circle (0.06/0.6);
              \end{scope}

             \draw[dashed] plot  (2-1/16,2/4+1/4) -- (2+1/16,2/4-1/4);
             \draw[dashed] plot  (1-1/16,1/4+1/4) -- (1+1/16,1/4-1/4);

       \draw plot (3,6/8) -- (6,0) ;
   
       \draw[dashed] plot  (5-1/16-1,2/4-1/4) -- (5+1/16-1,2/4+1/4);
       \draw[dashed] plot  (6-1/16-1,1/4-1/4) -- (6+1/16-1,1/4+1/4);

        \fill [black]   (3,6/8) circle (0.06);%
        
             \begin{scope}[shift={(0,-2)}]
              \draw plot (0,0) -- (3.5,7/8);
              \draw[dashed] plot  (2-1/16,2/4+1/4) -- (2+1/16,2/4-1/4);
              \draw[dashed] plot  (1-1/16,1/4+1/4) -- (1+1/16,1/4-1/4);
      
        \draw plot (2.5,7/8) -- (6,0);
       \draw[dashed] plot  (4-1/16,2/4-1/4) -- (4+1/16,2/4+1/4);
       \draw[dashed] plot  (5-1/16,1/4-1/4) -- (5+1/16,1/4+1/4);
              \draw   (3,6/8) circle (0.12);
        \fill [black]   (3,6/8) circle (0.06);%
             \end{scope}
            
             \end{tikzpicture}\end{center} \\ 
             \hline
       \end{tabular}
}

\def\FigB{
\begin{tabular}{|m{0.49\textwidth}|m{0.49\textwidth}|}
       \hline 
\begin{center}
\begin{tikzpicture}
       \draw plot (0,0) -- (4.5,9/10);
       \draw[dashed] plot  (2-1/16,2/5+1/4) -- (2+1/16,2/5-1/4);
       \draw[dashed] plot  (1-1/16,1/5+1/4) -- (1+1/16,1/5-1/4);
       \draw[dashed] plot  (3-1/16,3/5+1/4) -- (3+1/16,3/5-1/4);
       \draw[dashed] plot  (7/2-1/16,7/10+1/4) -- (7/2+1/16,7/10-1/4);
       \draw[dashed] plot  (5/2-1/16,5/10+1/4) -- (5/2+1/16,5/10-1/4);
       \draw[dashed] plot  (1.5-1/16,3/10+1/4) -- (1.5+1/16,3/10-1/4);
       \draw (4.5,9/10) node [right, scale=1] {2};
       \begin{scope}[shift={(0,-2)}]
              \draw plot (0,0) -- (4.5,9/10);
              \draw[dashed] plot  (2-1/16,2/5+1/4) -- (2+1/16,2/5-1/4);
              \draw[dashed] plot  (1-1/16,1/5+1/4) -- (1+1/16,1/5-1/4);
              \draw[dashed] plot  (3-1/16,3/5+1/4) -- (3+1/16,3/5-1/4);
              \draw[dashed] plot  (7/2-1/16,7/10+1/4) -- (7/2+1/16,7/10-1/4);
              \draw[dashed] plot  (5/2-1/16,5/10+1/4) -- (5/2+1/16,5/10-1/4);
              \draw[dashed] plot  (1.5-1/16,3/10+1/4) -- (1.5+1/16,3/10-1/4);
       \end{scope}
\end{tikzpicture}
\end{center}
&
\begin{center}
\begin{tikzpicture}

      \draw plot (1,1/5) -- (4.5,9/10);
      \draw[dashed] plot  (2-1/20,2/5+1/4) -- (2+1/20,2/5-1/4);
       \draw[dashed] plot  (2.5-1/20,5/10+1/4) -- (2.5+1/20,5/10-1/4);
       \draw[dashed] plot  (3-1/20,3/5+1/4) -- (3+1/20,3/5-1/4);
       \draw (1,1/5) node [left, scale=1] {1};
       \begin{scope}[shift={(8,0)}]
       \begin{scope}[yscale=1, xscale=-1]
              \draw plot (1,1/5) -- (4.5,9/10);
      \draw[dashed] plot  (2-1/20,2/5+1/4) -- (2+1/20,2/5-1/4);
       \draw[dashed] plot  (2.5-1/20,5/10+1/4) -- (2.5+1/20,5/10-1/4);
       \draw[dashed] plot  (3-1/20,3/5+1/4) -- (3+1/20,3/5-1/4);
       \draw (1,1/5) node [right, scale=1] {1};
       \end{scope}
\end{scope}
 \fill [black]   (4,6/8+1/32+1/64) circle (0.06);%

\begin{scope}[shift={(0,-2)}]
      \draw plot (1,1/5) -- (4.5,9/10);
      \draw[dashed] plot  (2-1/20,2/5+1/4) -- (2+1/20,2/5-1/4);
       \draw[dashed] plot  (2.5-1/20,5/10+1/4) -- (2.5+1/20,5/10-1/4);
       \draw[dashed] plot  (3-1/20,3/5+1/4) -- (3+1/20,3/5-1/4);
       \begin{scope}[shift={(8,0)}]
       \begin{scope}[yscale=1, xscale=-1]
              \draw plot (1,1/5) -- (4.5,9/10);
      \draw[dashed] plot  (2-1/20,2/5+1/4) -- (2+1/20,2/5-1/4);
       \draw[dashed] plot  (2.5-1/20,5/10+1/4) -- (2.5+1/20,5/10-1/4);
       \draw[dashed] plot  (3-1/20,3/5+1/4) -- (3+1/20,3/5-1/4);
       \end{scope}
\end{scope}
 \fill [black]   (4,6/8+1/32+1/64) circle (0.06);%
       
\end{scope}

\end{tikzpicture} \end{center}\\
\hline 
       \begin{center}
\begin{tikzpicture}

\draw plot (3/2,3/10) -- (5,1);
\draw (3/2,3/10) node [left, scale=1] {1};
              \draw[dashed] plot  (4-1/20-3/2,4/5+1/4-3/10) -- (4+1/20-3/2,4/5-1/4-3/10);
              \draw[dashed] plot  (4-1/20-2/2,4/5+1/4-2/10) -- (4+1/20-2/2,4/5-1/4-2/10);
              \draw[dashed] plot  (4-1/20-1/2,4/5+1/4-1/10) -- (4+1/20-1/2,4/5-1/4-1/10);
              \draw[dashed] plot  (4-1/20,4/5+1/4) -- (4+1/20,4/5-1/4);
              \begin{scope}[shift={(5,1)}, scale=0.6]
                     \draw plot[smooth, tension=1.3] coordinates {(0, 0) (0.5, 1/2) (-0.5, 1)};
                     \draw plot[smooth, tension=1.3] coordinates {(0, 0) (-0.5, 1/2) (0.5, 1)};
                     \fill [black]   (0,7/8) circle (0.06/0.6);
             \end{scope}
                   \fill [black]   (5,1) circle (0.06);%
             
              \draw plot (8-9/10,0+3/10) -- (5,1);
              \begin{scope}[shift={(-1/4, 1/12)}];
              \draw[dashed] plot  (7-1/12-1/2,1/3-1/4+1/6) -- (7+1/12-1/2,1/3+1/4+1/6);
              \draw[dashed] plot  (6-1/12,2/3-1/4) -- (6+1/12,2/3+1/4);
              \end{scope}

       \begin{scope}[shift={(0,-2)}]
              \draw plot (3/2,3/10) -- (5.5,1.1);
                            \draw[dashed] plot  (4-1/20-3/2,4/5+1/4-3/10) -- (4+1/20-3/2,4/5-1/4-3/10);
                            \draw[dashed] plot  (4-1/20-2/2,4/5+1/4-2/10) -- (4+1/20-2/2,4/5-1/4-2/10);
                            \draw[dashed] plot  (4-1/20-1/2,4/5+1/4-1/10) -- (4+1/20-1/2,4/5-1/4-1/10);
                            \draw[dashed] plot  (4-1/20,4/5+1/4) -- (4+1/20,4/5-1/4);

                     \draw plot (8-9/10,0+3/10) -- (4.5,1+1/6);
                     \begin{scope}[shift={(-1/4, 1/12)}];
                     \draw[dashed] plot  (7-1/12-1/2,1/3-1/4+1/6) -- (7+1/12-1/2,1/3+1/4+1/6);
                     \draw[dashed] plot  (6-1/12,2/3-1/4) -- (6+1/12,2/3+1/4);
                     \end{scope}

              \draw   (5,1) circle (0.12);
               \fill [black]   (5,1) circle (0.06);%
       \end{scope}
          
\end{tikzpicture}\end{center}
&
\begin{center}
\begin{tikzpicture}

 \draw (0,0) node [left, scale=1] {1};
 \draw plot (0,0) -- (3.5,7/8);
 \draw[dashed] plot  (2-1/16,2/4+1/4) -- (2+1/16,2/4-1/4);
 \draw[dashed] plot  (1-1/16,1/4+1/4) -- (1+1/16,1/4-1/4);
 \draw[dashed] plot  (3/2-1/16,3/8+1/4) -- (3/2+1/16,3/8-1/4);
  \fill [black]   (3,3/4) circle (0.06);%
 
   \draw plot (3,3/4) -- (4.5,1/8);
   \draw plot (3,3/4) -- (2.5,3/4+5/24);
\draw[dashed] plot  (3+3/4+5/12*1/4,3/4-5/12*3/4+1/4) -- (3+3/4-5/12*1/4,3/4-5/12*3/4-1/4);
   \begin{scope}[shift={(4,0)}]
         \draw plot (0.5,0.125) -- (3,6/8);
         \draw[dashed] plot  (3/2+1/2-1/16,2/4+1/4) -- (3/2+1/16+1/2,2/4-1/4);
         \draw[dashed] plot  (1-1/16+1/2,1/4+1/4+1/8) -- (1+1/16+1/2,1/4-1/4+1/8);
  \fill [black]   (1/2,1/8) circle (0.06);%
 \begin{scope}[shift={(1/2,1/8)}, scale=0.6,yscale=-1]
       \draw plot[smooth, tension=1.5] coordinates {(0, 0) (0.5, 1/2) (-0.5, 1)};
       \draw plot[smooth, tension=1.5] coordinates {(0, 0) (-0.5, 1/2) (0.5, 1)};
       \fill [black]   (0,7/8) circle (0.06/0.6);
\end{scope}
   \end{scope}

      \begin{scope}[shift={(0,-2)}]
      \draw plot (0,0) -- (3.5,7/8);
      \draw[dashed] plot  (2-1/16,2/4+1/4) -- (2+1/16,2/4-1/4);
      \draw[dashed] plot  (1-1/16,1/4+1/4) -- (1+1/16,1/4-1/4);
      \draw[dashed] plot  (3/2-1/16,3/8+1/4) -- (3/2+1/16,3/8-1/4);
       \fill [black]   (3,3/4) circle (0.06);%
      
        \draw plot (3,3/4) -- (4.5,1/8);
        \draw plot (5,1/8-5/24) -- (4.5,1/8);
        \draw plot (3,3/4) -- (2.5,3/4+5/24);
     \draw[dashed] plot  (3+3/4+5/12*1/4,3/4-5/12*3/4+1/4) -- (3+3/4-5/12*1/4,3/4-5/12*3/4-1/4);
        \begin{scope}[shift={(4,0)}]
              \draw plot (0,0) -- (3,6/8);
              \draw[dashed] plot  (3/2+1/2-1/16,2/4+1/4) -- (3/2+1/16+1/2,2/4-1/4);
              \draw[dashed] plot  (1-1/16+1/2,1/4+1/4+1/8) -- (1+1/16+1/2,1/4-1/4+1/8);
       \fill [black]   (1/2,1/8) circle (0.06);%
       
       \draw   (1/2,1/8) circle (0.12);
        \end{scope}
      \end{scope}
    
\end{tikzpicture}
\end{center}
 \\
\hline
\begin{center}
\begin{tikzpicture}
       
              \draw plot (0.5,1/8) -- (3,6/8);
              \draw[dashed] plot  (2-1/16,2/4+1/4) -- (2+1/16,2/4-1/4);
              \draw[dashed] plot  (1-1/16+1/2,1/4+1/4+1/8) -- (1+1/16+1/2,1/4-1/4+1/8);
               \fill [black]   (3,3/4) circle (0.06);%
              \begin{scope}[shift={(3,3/4)}, scale=0.6]
                      \draw plot[smooth, tension=1.3] coordinates {(0, 0) (0.5, 1/2) (-0.5, 1)};
                      \draw plot[smooth, tension=1.3] coordinates {(0, 0) (-0.5, 1/2) (0.5, 1)};
                      \fill [black]   (0,7/8) circle (0.06/0.6);
              \end{scope}

                \draw plot (3,3/4) -- (5.5,1/8);
              
                \draw[dashed] plot  (5-1/2-1/16,1/4-1/4+1/8) -- (5+1/16-1/2,1/4+1/4+1/8);
                \draw[dashed] plot  (4-1/16,2/4-1/4) -- (4+1/16,2/4+1/4);
                \begin{scope}[shift={(5,0)}]
                      \draw plot (1/2,1/8) -- (3,6/8);
                      \draw[dashed] plot  (2+1/2-1/16-1/2,2/4+1/4+1/8-1/8) -- (2+1/16+1/2-1/2,2/4-1/4+1/8-1/8);
                      \draw[dashed] plot  (1-1/16+1/2,1/4+1/4+1/8) -- (1+1/16+1/2,1/4-1/4+1/8);
               \fill [black]   (1/2,1/8) circle (0.06);%
              \begin{scope}[shift={(1/2,1/8)}, scale=0.6,yscale=-1]
                     \draw plot[smooth, tension=1.3] coordinates {(0, 0) (0.5, 1/2) (-0.5, 1)};
                     \draw plot[smooth, tension=1.3] coordinates {(0, 0) (-0.5, 1/2) (0.5, 1)};
                     \fill [black]   (0,7/8) circle (0.06/0.6);
             \end{scope}
             
                \end{scope}
       
             \begin{scope}[shift={(0,-2)}]
             \draw plot (0.5,1/8) -- (3.5,7/8);
             \draw[dashed] plot  (2-1/16,2/4+1/4) -- (2+1/16,2/4-1/4);
             \draw[dashed] plot  (1-1/16+1/2,1/4+1/4+1/8) -- (1+1/16+1/2,1/4-1/4+1/8);
              \fill [black]   (3,3/4) circle (0.06);%
               \draw   (3,3/4) circle (0.12);
               \draw plot (3,3/4) -- (6,0);
               \draw plot (3,3/4) -- (2.5,7/8);
               \draw[dashed] plot  (5-1/2-1/16,1/4-1/4+1/8) -- (5+1/16-1/2,1/4+1/4+1/8);
               \draw[dashed] plot  (4-1/16,2/4-1/4) -- (4+1/16,2/4+1/4);
               \begin{scope}[shift={(5,0)}]
                     \draw plot (0,0) -- (3,6/8);
                     \draw[dashed] plot  (2+1/2-1/16-1/2,2/4+1/4+1/8-1/8) -- (2+1/16+1/2-1/2,2/4-1/4+1/8-1/8);
                     \draw[dashed] plot  (1-1/16+1/2,1/4+1/4+1/8) -- (1+1/16+1/2,1/4-1/4+1/8);
              \fill [black]   (1/2,1/8) circle (0.06);%
              \draw   (1/2,1/8) circle (0.12);
               \end{scope}
             \end{scope}
           
       \end{tikzpicture}
\end{center}
&
\begin{center}
\begin{tikzpicture}
\draw plot (1/2,1/8) -- (3,6/8);
\draw[dashed] plot  (2-1/16,2/4+1/4) -- (2+1/16,2/4-1/4);
\draw[dashed] plot  (1+1/2-1/16,1/4+1/4+1/8) -- (1+1/2+1/16,1/4-1/4+1/8);
 \fill [black]   (3,6/8) circle (0.06);%
\begin{scope}[shift={(3,3/4)}, scale=0.6]
       \draw plot[smooth, tension=1.3] coordinates {(0, 0) (0.5, 1/2) (-0.5, 1)};
       \draw plot[smooth, tension=1.3] coordinates {(0, 0) (-0.5, 1/2) (0.5, 1)};
       \fill [black]   (0,7/8) circle (0.06/0.6);
\end{scope}

\draw plot (9/2,0) -- (6/2,6/8);
\draw[dashed] plot  (7/2+1/8+1/8,-1/16+4/8+1/4) -- (7/2+1/8-1/8,4/8-1/4-1/16);
 \fill [black]   (8/2+1/4,1/8) circle (0.06);%

\begin{scope}[yscale=1, xscale=-1, shift={(-17/2,0)}]
        \draw plot (1/2,1/8) -- (3,6/8);
        \draw[dashed] plot  (2-1/16,2/4+1/4) -- (2+1/16,2/4-1/4);
        \draw[dashed] plot  (1+1/2-1/16,1/4+1/4+1/8) -- (1+1/2+1/16,1/4-1/4+1/8);
 \fill [black]   (3,6/8) circle (0.06);%

 \begin{scope}[shift={(3,3/4)}, scale=0.6]
       \draw plot[smooth, tension=1.3] coordinates {(0, 0) (0.5, 1/2) (-0.5, 1)};
       \draw plot[smooth, tension=1.3] coordinates {(0, 0) (-0.5, 1/2) (0.5, 1)};
       \fill [black]   (0,7/8) circle (0.06/0.6);
\end{scope}

\draw plot (9/2,0) -- (6/2,6/8);
\draw[dashed] plot  (7/2+1/8+1/8,-1/16+4/8+1/4) -- (7/2+1/8-1/8,4/8-1/4-1/16);
\end{scope}

      \begin{scope}[shift={(0,-2)}]
\draw plot (7/2,7/8) -- (3,6/8);
\draw plot (1/2,1/8) -- (3,6/8);
\draw[dashed] plot  (2-1/16,2/4+1/4) -- (2+1/16,2/4-1/4);
\draw[dashed] plot  (1+1/2-1/16,1/4+1/4+1/8) -- (1+1/2+1/16,1/4-1/4+1/8);
 \fill [black]   (3,6/8) circle (0.06);%
 \draw   (3,6/8) circle (0.12);
\draw plot (9/2,0) -- (6/2,6/8);
\draw plot (5/2,8/8) -- (6/2,6/8);
\draw[dashed] plot  (7/2+1/8+1/8,-1/16+4/8+1/4) -- (7/2+1/8-1/8,4/8-1/4-1/16);
 \fill [black]   (8/2+1/4,1/8) circle (0.06);%

\begin{scope}[yscale=1, xscale=-1, shift={(-17/2,0)}]
        \draw plot (7/2,7/8) -- (3,6/8);
        \draw plot (1/2,1/8) -- (3,6/8);
        \draw[dashed] plot  (2-1/16,2/4+1/4) -- (2+1/16,2/4-1/4);
        \draw[dashed] plot  (1+1/2-1/16,1/4+1/4+1/8) -- (1+1/2+1/16,1/4-1/4+1/8);
 \fill [black]   (3,6/8) circle (0.06);%
 \draw   (3,6/8) circle (0.12);
\draw plot (9/2,0) -- (6/2,6/8);
\draw plot (5/2,8/8) -- (6/2,6/8);
\draw[dashed] plot  (7/2+1/8+1/8,-1/16+4/8+1/4) -- (7/2+1/8-1/8,4/8-1/4-1/16);
\end{scope}
      \end{scope}
\end{tikzpicture} \end{center}\\
\hline
\begin{center}
\begin{tikzpicture}
      \draw plot (0,0) -- (2*3,0);
      \draw plot (0,0.5) -- (2*3,0.5);
       \draw plot (2*0.5, 1) -- (2*0.5, -1);
       \draw[dashed] plot  (2*0.5-1/4, -0.5) -- (2*0.5+1/4, -0.5);
       \draw[dashed] plot  (2*0.5-1/4, -3/4) -- (2*0.5+1/4, -3/4);
       \fill [black]   (2*1/2,0) circle (0.06);%
       \fill [black]   (2*1/2,1/2) circle (0.06);%
      
      \begin{scope}[shift={(2*1,0)}]
      \draw plot (2*0.5, 1) -- (2*0.5, -1);
      \draw[dashed] plot  (2*0.5-1/4, -0.5) -- (2*0.5+1/4, -0.5);
      \draw[dashed] plot  (2*0.5-1/4, -3/4) -- (2*0.5+1/4, -3/4);
      \fill [black]   (2*1/2,0) circle (0.06);%
      \fill [black]   (2*1/2,1/2) circle (0.06);%
    
      \end{scope}
        \begin{scope}[shift={(2*2,0)}]
              \draw plot (2*0.5, 1) -- (2*0.5, -1);
              \draw[dashed] plot  (2*0.5-1/4, -0.5) -- (2*0.5+1/4, -0.5);
              \draw[dashed] plot  (2*0.5-1/4, -3/4) -- (2*0.5+1/4, -3/4);
              \fill [black]   (2*1/2,0) circle (0.06);%
              \fill [black]   (2*1/2,1/2) circle (0.06);%
    
             \end{scope}

       \begin{scope}[shift={(0,-2)}]
      \draw plot (0,0) -- (2*3,0);
       \draw plot (2*0.5, 0.5) -- (2*0.5, -1);
       \draw[dashed] plot  (2*0.5-1/4, -0.5) -- (2*0.5+1/4, -0.5);
       \draw[dashed] plot  (2*0.5-1/4, -3/4) -- (2*0.5+1/4, -3/4);
       \fill [black]   (2*1/2,0) circle (0.06);%
       \draw   (2*1/2,0) circle (0.12);
       \begin{scope}[shift={(2*1,0)}]
              \draw plot (2*0.5, 0.5) -- (2*0.5, -1);
              \draw[dashed] plot  (2*0.5-1/4, -0.5) -- (2*0.5+1/4, -0.5);
              \draw[dashed] plot  (2*0.5-1/4, -3/4) -- (2*0.5+1/4, -3/4);
              \fill [black]   (2*1/2,0) circle (0.06);%
              \draw   (2*1/2,0) circle (0.12);
       \end{scope}
         \begin{scope}[shift={(2*2,0)}]
              \draw plot (2*0.5, 0.5) -- (2*0.5, -1);
              \draw[dashed] plot  (2*0.5-1/4, -0.5) -- (2*0.5+1/4, -0.5);
              \draw[dashed] plot  (2*0.5-1/4, -3/4) -- (2*0.5+1/4, -3/4);
              \fill [black]   (2*1/2,0) circle (0.06);%
              \draw   (2*1/2,0) circle (0.12);
              \end{scope}

       \end{scope}
\end{tikzpicture} \end{center}
&
  \\
\hline
\end{tabular}
}

\def\FigC{
\begin{tikzpicture}[
              outer frame sep=0ex,%
             /tikz/inner frame xsep=4ex,
             /tikz/inner frame ysep=2ex,
             show background top,%
             show background bottom,%
             show background left,%
             show background right]
             \draw plot (0,0) -- (3,0);
             \draw[domain=0:1] plot ({\x+3)},{\x^2/2});
             \begin{scope}[shift={(0, -3/4)}]
              \draw plot (0,0) -- (3,0);
              \draw[domain=0:1] plot ({\x+3)},{\x^2/2});
             \end{scope}
             \draw plot (1,-1) -- (1, 1.5);
             \draw[dashed] plot  (2*0.5-1/4, 0.5) -- (2*0.5+1/4, 0.5);
             \draw[dashed] plot  (2*0.5-1/4,1) -- (2*0.5+1/4, 1);
             \draw plot (2,-1) -- (2, 1.5);
             \draw[dashed] plot  (1+2*0.5-1/4, 0.5) -- (1+2*0.5+1/4, 0.5);
             \draw[dashed] plot  (1+2*0.5-1/4,1) -- (1+2*0.5+1/4, 1);
             \fill [black]   (1,0) circle (0.06);%
             \fill [black]   (1,-3/4) circle (0.06);%
             \fill [black]   (2,-3/4) circle (0.06);%
             \fill [black]   (2,0) circle (0.06);%
             \fill [black]   (7/2,1/8-3/4) circle (0.06);%
             \fill [black]   (7/2,1/8) circle (0.06);%
             \begin{scope}[xscale=-1, shift={(-7,0)}]
             \draw plot (0,0) -- (3,0);
             \draw[domain=0:1] plot ({\x+3)},{\x^2/2});
             \begin{scope}[shift={(0, -3/4)}]
              \draw plot (0,0) -- (3,0);
              \draw[domain=0:1] plot ({\x+3)},{\x^2/2});
             \end{scope}
             \draw plot (1,-1) -- (1, 1.5);
             \draw[dashed] plot  (2*0.5-1/4, 0.5) -- (2*0.5+1/4, 0.5);
             \draw[dashed] plot  (2*0.5-1/4,1) -- (2*0.5+1/4, 1);
             \draw plot (2,-1) -- (2, 1.5);
             \draw[dashed] plot  (1+2*0.5-1/4, 0.5) -- (1+2*0.5+1/4, 0.5);
             \draw[dashed] plot  (1+2*0.5-1/4,1) -- (1+2*0.5+1/4, 1);
             \fill [black]   (1,0) circle (0.06);%
             \fill [black]   (1,-3/4) circle (0.06);%
             \fill [black]   (2,-3/4) circle (0.06);%
             \fill [black]   (2,0) circle (0.06);%
             \end{scope}

             \begin{scope}[shift={(0, -3)}]
             \draw plot (0,0) -- (3,0);
              \draw[domain=0:1] plot ({\x+3)},{\x^2/2});
              \draw plot (1,-0.5) -- (1, 1.5);
              \draw[dashed] plot  (2*0.5-1/4, 0.5) -- (2*0.5+1/4, 0.5);
              \draw[dashed] plot  (2*0.5-1/4,1) -- (2*0.5+1/4, 1);
              \draw plot (2,-0.5) -- (2, 1.5);
              \draw[dashed] plot  (1+2*0.5-1/4, 0.5) -- (1+2*0.5+1/4, 0.5);
              \draw[dashed] plot  (1+2*0.5-1/4,1) -- (1+2*0.5+1/4, 1);
              \fill [black]   (1,0) circle (0.06);%
              \draw   (2*1/2,0) circle (0.12);
              \fill [black]   (2,0) circle (0.06);%
              \draw   (2,0) circle (0.12);
              \fill [black]   (7/2,1/8) circle (0.06);%
              \draw   (7/2,1/8) circle (0.12);
              \begin{scope}[xscale=-1, shift={(-7,0)}]
                     \draw plot (0,0) -- (3,0);
                     \draw[domain=0:1] plot ({\x+3)},{\x^2/2});
              \draw plot (1,-0.5) -- (1, 1.5);
              \draw[dashed] plot  (2*0.5-1/4, 0.5) -- (2*0.5+1/4, 0.5);
              \draw[dashed] plot  (2*0.5-1/4,1) -- (2*0.5+1/4, 1);
              \draw plot (2,-0.5) -- (2, 1.5);
              \draw[dashed] plot  (1+2*0.5-1/4, 0.5) -- (1+2*0.5+1/4, 0.5);
              \draw[dashed] plot  (1+2*0.5-1/4,1) -- (1+2*0.5+1/4, 1);
              \fill [black]   (1,0) circle (0.06);%
              \draw   (2*1/2,0) circle (0.12);
              \fill [black]   (2,0) circle (0.06);%
              \draw   (2,0) circle (0.12);
              \end{scope}

              \end{scope}
\end{tikzpicture}
}

\def\FigD{
\begin{tikzpicture}[
       outer frame sep=0ex,%
      /tikz/inner frame xsep=4ex,
      /tikz/inner frame ysep=2ex,
      show background top,%
      show background bottom,%
      show background left,%
      show background right]
\draw plot (0,0) -- (3,0);
\draw plot (0,0.5) -- (3,0.5);
\draw plot (1,-0.5) -- (1, 3);
\draw[dashed] plot  (2*0.5-1/4, 2.5) -- (2*0.5+1/4, 2.5);
\draw[dashed] plot  (2*0.5-1/4,1) -- (2*0.5+1/4, 1);
\draw[dashed] plot  (2*0.5-1/4,1.5) -- (2*0.5+1/4, 1.5);
\draw[dashed] plot  (2*0.5-1/4,2) -- (2*0.5+1/4, 2);
\draw (1,3) node [left, scale=1] {1};
\draw plot (2,-0.5) -- (2, 2);
\draw[dashed] plot  (1+2*0.5-1/4, 1.5) -- (1+2*0.5+1/4, 1.5);
\draw[dashed] plot  (1+2*0.5-1/4,1) -- (1+2*0.5+1/4, 1);
\draw[domain=pi/2:3*pi/2] plot ({1/2*sin(\x r)+3.25},{1/2*cos(\x r)});
\draw plot (2.75,0) -- (2.75, 0.75);
\draw plot (3.75,0) -- (3.75, 0.25);
\draw[dashed] plot  (3.25, -0.25) -- (3.25, -0.75);
\draw plot (3.5,0) -- (5,0);
              \draw[domain=pi/2:pi] plot ({1/2*sin(\x r)+5},{1/2*cos(\x r)+1/2});
              \draw plot (5.5,1/2) -- (5.5,3/2);
\draw plot (4.5,-0.5) -- (4.5, 2);
\draw (4.5,2) node [left, scale=1] {1};
\draw[dashed] plot  (4.5-1/4, 0.5) -- (4.5+1/4, 0.5);
\draw[dashed] plot  (4.5-1/4,1) -- (4.5+1/4, 1);
\draw[dashed] plot  (4.5-1/4,1.5) -- (4.5+1/4, 1.5);
\draw plot (5.5,3/4) -- (8.5, 3/4);
\draw (8.5, 3/4) node [right, scale=1] {1};
\draw[domain=pi:4/2*pi] plot ({1/4*sin(\x r)+5.5},{1/4*cos(\x r)+1});
\draw plot (5.5, 5/4) -- (5.75, 5/4);
\draw[dashed] plot  (6,3/4+1/4) -- (6, 3/4-1/4);
\draw[dashed] plot  (6.5,3/4+1/4) -- (6.5, 3/4-1/4);
\draw[dashed] plot  (7,3/4+1/4) -- (7, 3/4-1/4);
\draw plot (7.75,2+1/4) -- (7.75, -0.25);
\draw (7.75,2+1/4) node [left, scale=1] {1};
\draw[dashed] plot  (7.75-1/4,1.5+1/4) -- (7.75+1/4, 1.5+1/4);
\draw[dashed] plot  (7.75-1/4,2-1+0.25) -- (7.75+1/4, 2-1+1/4);
\draw[dashed] plot  (7.75-1/4,+0.25) -- (7.75+1/4, +1/4);

 \fill [black]   (1,0) circle (0.06);%
 \fill [black]   (1,0.5) circle (0.06);%
 \fill [black]   (2,0) circle (0.06);%
 \fill [black]   (2,0.5) circle (0.06);%
 \fill [black]   (2.75,0) circle (0.06);%
 \fill [black]   (2.75,0.5) circle (0.06);%
 \fill [black]   (3.75,0) circle (0.06);%
 \fill [black]   (4.5,0) circle (0.06);%
 \fill [black]   (5.5,3/4) circle (0.06);%
 \fill [black]   (5.5, 5/4) circle (0.06);%
 \fill [black]   (7.75,3/4) circle (0.06);%
 
\begin{scope}[shift={(0,-3.5 )}]
\draw plot (0,0) -- (3,0);
\draw plot (1,-0.5) -- (1, 2.5);
\draw[dashed] plot  (2*0.5-1/4, 0.5) -- (2*0.5+1/4, 0.5);
\draw[dashed] plot  (2*0.5-1/4,1) -- (2*0.5+1/4, 1);
\draw[dashed] plot  (2*0.5-1/4,1.5) -- (2*0.5+1/4, 1.5);
\draw[dashed] plot  (2*0.5-1/4,2) -- (2*0.5+1/4, 2);
\draw plot (2,-0.5) -- (2, 1.5);
\draw[dashed] plot  (1+2*0.5-1/4, 0.5) -- (1+2*0.5+1/4, 0.5);
\draw[dashed] plot  (1+2*0.5-1/4,1) -- (1+2*0.5+1/4, 1);
\draw[domain=pi/2:3*pi/2] plot ({1/2*sin(\x r)+3.25},{1/2*cos(\x r)});
\draw plot (2.75,0) -- (2.75, 0.25);
\draw plot (3.75,0) -- (3.75, 0.25);
\draw[dashed] plot  (3.25, -0.25) -- (3.25, -0.75);
\draw plot (3.5,0) -- (5,0);
              \draw[domain=pi/2:pi] plot ({1/2*sin(\x r)+5},{1/2*cos(\x r)+1/2});
              \draw plot (5.5,1/2) -- (5.5,2/2);
\draw plot (4.5,-0.5) -- (4.5, 2);
\draw[dashed] plot  (4.5-1/4, 0.5) -- (4.5+1/4, 0.5);
\draw[dashed] plot  (4.5-1/4,1) -- (4.5+1/4, 1);
\draw[dashed] plot  (4.5-1/4,1.5) -- (4.5+1/4, 1.5);
\draw plot (5.25,3/4) -- (8.5, 3/4);
\draw[dashed] plot  (6,3/4+1/4) -- (6, 3/4-1/4);
\draw[dashed] plot  (6.5,3/4+1/4) -- (6.5, 3/4-1/4);
\draw[dashed] plot  (7,3/4+1/4) -- (7, 3/4-1/4);
\draw plot (7.75,2+1/4) -- (7.75, -0.25);
\draw[dashed] plot  (7.75-1/4,1.5+1/4) -- (7.75+1/4, 1.5+1/4);
\draw[dashed] plot  (7.75-1/4,2-1+0.25) -- (7.75+1/4, 2-1+1/4);
\draw[dashed] plot  (7.75-1/4,+0.25) -- (7.75+1/4, +1/4);

 \fill [black]   (1,0) circle (0.06);%
 \draw   (2*1/2,0) circle (0.12);
 \fill [black]   (2,0) circle (0.06);%
 \draw   (2,0) circle (0.12);
 \fill [black]   (2.75,0) circle (0.06);%
 \draw   (2.75,0) circle (0.12);
 \fill [black]   (3.75,0) circle (0.06);%
 \fill [black]   (4.5,0) circle (0.06);%
 \fill [black]   (5.5,3/4) circle (0.06);%
 \draw   (5.5,3/4) circle (0.12);
 \fill [black]   (7.75,3/4) circle (0.06);%
\end{scope}
\end{tikzpicture}
}